\documentclass[autocontact]{gaceta}

\usepackage[spanish]{babel}

\usepackage[applemac]{inputenc}

\theoremstyle{defnition}

\newtheorem{teo}{\bf TEOREMA}
\newtheorem{cor}[teo]{\bf COROLARIO}
\newtheorem{concepto}{\bf DEFINICI\'{O}N}
\newtheorem{lem}[teo]{\bf LEMA}

\newtheorem{remark}[teo]{NOTA}

\journame{La Gaceta de la RSME}
\yearofpublication{0000}
\volume{00}
\issuenumber{0}

\belongstopart{Art\'{\i}culos} 

\title{Los teoremas de Fr\'{e}chet,  Montel y Popoviciu y los grafos de los polinomios discontinuos}
\author{J. M. Almira y Kh. F. Abu-Helaiel} 
\shorttitle{ teoremas de Fr\'{e}chet,  Montel y Popoviciu}

\contact{J. M. Almira, Dpto. de Matem\'aticas, Universidad de Jan}
{jmalmira@ujaen.es}{}
\contact{Kh. F. Abu-Helaiel, Dpto. de Matem\'aticas, Universidad de Jan}
{kabu@ujaen.es}{}

\begin{document}

\maketitle

\begin{abstract}
Pretendemos introducir al lector en los problemas de regularidad para las ecuaciones funcionales. Para ello, hacemos un seguimiento detallado de algunos resultados relacionados con la regularidad de la ecuacin funcional de Frchet. En particular, se demuestran los teoremas clsicos de Frchet, Montel y Montel-Popoviciu, y se estudia el grafo de las soluciones discontinuas de la ecuacin funcional de Frchet. Este artculo tiene vocacin de ser un homenaje a Tiberiu Popoviciu, un matemtico rumano que, en nuestra opinin, merece ser recordado entre los grandes del anlisis matemtico. 

\end{abstract}

\section{Motivacin: la ecuacin funcional de Cauchy}

Hace algunos aos apareci, en esta Gaceta de la RSME,  un interesante artculo sobre ecuaciones funcionales  \cite{castillo}. En dicho trabajo el profesor Enrique Castillo, de la Universidad de Cantabria, nos 
convenci a todos los lectores (o, cuando menos, al primero de los autores de este trabajo, quien descubrira esta joya poco tiempo despus,  cuando particip en la 
Redaccin de esta revista), del inters fundamental 
que tiene esta rama del anlisis matemtico, no slo por su belleza sino tambin por su enorme versatilidad. Tanto fue as, que algunos nos pusimos manos a la obra y a partir de entonces hemos dedicado una parte esencial de nuestros esfuerzos en investigar sobre estos temas. En aquel artculo se enfatiz mucho el papel de los modelos. Ahora, con esta nueva contribucin, nos gustara presentar un rpido repaso del estado del arte
en relacin a una de las ecuaciones funcionales clsicas por excelencia: la ecuacin de Frchet. Con ello, pretendemos introducir al lector, mediante la exposicin detallada de una serie de resultados clsicos, as como algunas contribuciones recientes, a uno de los temas fundamentales de esta teora: los problemas de regularidad. En su artculo \cite{castillo} el autor nos recordaba que algunos de los nombres ms importantes de los siglos XVIII, XIX y principios del siglo XX, como son (en orden cronolgico) DÕAlembert, Euler, Gauss, Cauchy, Abel, Weierstrass, Darboux o Hilbert, trabajaron por algn tiempo con ecuaciones funcionales. Nosotros queremos ampliar esta lista al incluir ahora los nombres de Frchet, Montel y Popoviciu. 

\begin{figure}
  \centering
    \includegraphics[scale=0.49]{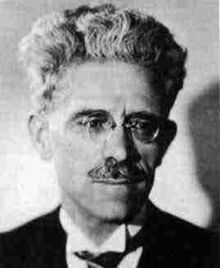} \includegraphics[scale=0.4]{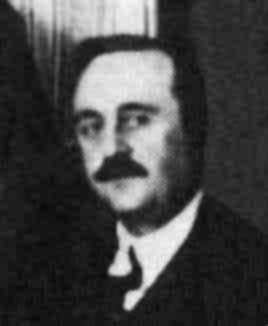} \includegraphics[scale=0.5]{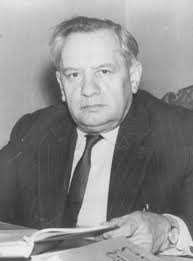} 
  \caption{Maurice Frchet, Paul Montel y Tiberiu Popoviciu}
\end{figure}

Como es natural, puesto que las ecuaciones funcionales no forman parte del curriculum que se estudia en los grados de matemticas, nos gustara comenzar explicando algunos de los resultados clsicos que motivaron la teora.  Ya hemos formulado nuestra opinin de esta parte del anlisis, pero ahora vamos a intentar mostrar que, en efecto, esta teora es extremadamente hermosa. Para ello, consideramos la ecuacin funcional ms clsica de todas: la ecuacin de Cauchy:  
\begin{equation}\label{cauchy}
f(x+y)=f(x)+f(y).
\end{equation} 
Las soluciones $f$ de la ecuacin  \eqref{cauchy} se llaman, por motivos evidentes, funciones aditivas. Cauchy demostr que si $f:\mathbb{R}\to \mathbb{R}$ es una funcin continua que satisface \eqref{cauchy}, entonces  $f(x)=ax$ para cierta constante $a$. Veamos cmo se puede demostrar esto. La aditividad de $f$ nos permite afirmar que, si $x\in\mathbb{R}$ y $n\in\mathbb{N}$, entonces $f(nx)=nf(x)$. Por tanto, $f(x)=f(n(x/n))=nf(x/n)$, y, en consecuencia, $f(x/n)=f(x)/n$ para todo $x$ real y todo nmero natural $n$. Adems, $f(0)=2f(0)$ implica que $f(0)=0$, y $f(x)+f(-x)=f(0)=0$ implica que $f(-x)=-f(x)$, por lo que las propiedades anteriores se trasladan inmediatamente al caso $x\in\mathbb{R}$, $n\in\mathbb{Z}$. Se sigue que, si $r=n/m\in\mathbb{Q}$, entonces $f(rx)=f(n(x/m))=nf(x/m)=(n/m)f(x)=rf(x)$ para todo nmero real $x$. Tomando ahora $x=1$, concluimos que $f(r)=rf(1)$ para todo $r\in\mathbb{Q}$. Como $\mathbb{Q}$ es denso en $\mathbb{R}$ y $f$ es continua, concluimos que $f(x)=f(1)x=ax$ con $a=f(1)$, para todo $x\in\mathbb{R}$.  

Algunos aos despus, Darboux \cite{darboux} (ver tambin \cite{AD})  demostrara que si una funcin aditiva es continua en un punto, entonces es continua en todos los puntos de la recta. Esto es consecuencia de la  aditividad de $f$ , pues si tenemos garantizada la continuidad en un punto $x_0$ y tomamos otro punto $x$ cualquiera, entonces, para todo $h\in\mathbb{R}$, 
\begin{eqnarray*}
|f(x+h)-f(x)| &=& |f(x+x_0+h)-f(x_0)-f(x)|\\
&=& |f(x)+f(x_0+h)-f(x_0)-f(x)| = |f(x_0+h)-f(x_0)|.
\end{eqnarray*}
Es ms, Darboux demostr que si $f$ est acotada en algn intervalo abierto no vaco, entonces $f(x)=f(1)x$ para todo $x\in\mathbb{R}$. Este mismo resultado, en una expresin an ms contundente, sera demostrado por Ricardo San Juan \cite{sanjuan}, uno de los primeros matemticos espaoles que alcanzaron el prestigio internacional en el siglo XX. Su argumento es  tan elegante que no podemos permitirnos aqu obviarlo. Adems, como veremos luego, motiv buena parte de nuestro trabajo en ecuaciones funcionales. El profesor San Juan demostr el siguiente resultado:

\begin{teo}[R. San Juan] Si $f:\mathbb{R}\to\mathbb{R}$ es una funcin aditiva discontinua, entonces su grafo, $G(f)=\{(x,f(x)):x\in\mathbb{R}\}$ es un subconjunto denso del plano.
\end{teo} 

\noindent \textbf{Demostracin. } Supongamos que  $f$ es aditiva pero no es una funcin del tipo  $f(x)=ax$ para ninguna constante  $a$ (ya sabemos que las soluciones continuas de la ecuacin de Cauchy son de este tipo). Entonces existen   $x_0,x_1\neq 0$ 
tales que  $f(x_0)/x_0 \neq f(x_1)/x_1$. En otras palabras, 
$f(x_0)x_1-f(x_1)x_0\neq 0$. Esto significa que el determinante de la matriz 
$$A=\left[\begin{array}{ccc} x_0 & x_1\\   f(x_0) & f(x_1)\end{array}\right]$$
es distinto de cero. Por tanto, los vectores 
$(x_0,f(x_0)), (x_1,f(x_1))$ 
son linealmente independientes. Se sigue que las combinaciones lineales que podemos realizar con estos vectores, tomando coeficientes en $\mathbb{Q}$, forman un subconjunto denso del plano. Ahora bien,  
si  $r_0,r_1 \in \mathbb{Q}$, entonces
\[
r_0(x_0,f(x_0))+r_1(x_1,f(x_1))=(r_0x_0+r_1x_1,f(r_0x_0+r_1x_1))\in G(f).
\]
Por tanto, $G(f)$ es denso en $\mathbb{R}^2$ {\hfill $\Box$}

En 1906, Hamel \cite{hamel} haba introducido sus bases (es decir, las bases algebraicas de $\mathbb{R}$ como espacio vectorial sobre $\mathbb{Q}$, cuya existencia est garantizada por el axioma de eleccin) precisamente para demostrar la existencia de soluciones discontinuas de la ecuacin de Cauchy.  En efecto, si $\beta=\{v_i\}_{i\in I}$ es una de tales bases y $\varphi:\beta\to \mathbb{R}$ es una aplicacin arbitraria, entonces la nica aplicacin $\mathbb{Q}$-lineal $\mathcal{L}_{\beta,\varphi}:\mathbb{R}\to\mathbb{R}$ que verifica $\mathcal{L}(v_i)=\varphi(i)$ para todo $i\in I$, es una funcin aditiva, y esta funcin ser discontinua si (y solo si) existen $i_0,i_1$ tales que $ \varphi(i_0)/v_{i_0}\neq  \varphi(i_1)/v_{i_1}$, cosa que podemos forzar sin problemas.
 
Posteriormente, Sierpinsky \cite{sierpinsky} y Banach \cite{banach} publicaron, en el primer volumen de la revista Fundamenta Mathematicae, dos demostraciones diferentes de que las funciones $f:\mathbb{R}\to\mathbb{R}$ aditivas medibles son necesariamente de la forma $f(x)=ax$ para cierta constante $a$. Por ltimo,  Kormes \cite{kormes} utiliz el teorema de Darboux y la propiedad de los conjuntos medibles de Lebesgue -que haba sido demostrada recientemente por Steinhauss \cite[Teorema VII]{steinhaus}- segn la cual, si $A\subseteq \mathbb{R}$  es un conjunto con medida de Lebesgue positiva, $|A|>0$, entonces el conjunto  
$$A+A=\{x+y:x,y \in A\}$$ contiene un intervalo abierto no vaco, para demostrar el siguiente resultado, ms fuerte que los obtenidos por Banach y Sierpinsky:

\begin{teo}[Kormes] \label{kormes} Si  $f$ es una funcin aditiva y acotada en un conjunto $A\subseteq \mathbb{R}$ con medida de Lebesgue positiva, $|A|>0$, entonces $f(x)=f(1)x$ para todo $x\in\mathbb{R}$.
\end{teo}

\noindent \textbf{Demostracin. } Si $\sup_{x\in A}|f(x)|\leq M$ entonces $\sup_{x\in A+A}|f(x)|\leq 2M$, pues $f$ es aditiva y, por tanto, $|f(x+y)|=|f(x)+f(y)|\leq |f(x)|+|f(y)|$ para todo $x,y\in A$. Se sigue que $f$ est acotada en un intervalo abierto no vaco y  el teorema de Darboux implica  que $f(x)=ax$. {\hfill $\Box$}

Como se ve, hasta ahora las cosas han salido fciles. La ecuacin de Cauchy ha demostrado ser muy productiva. Ella posee la sorprendente cualidad de que, sin aparecer ningn tipo de condiciones de regularidad  en su definicin (las funciones $f$ que satisfacen la ecuacin pueden ser altamente irregulares), resulta que si exigimos un mnimo de regularidad sobre una de sus soluciones (como, por ejemplo, estar acotada en un conjunto de medida positiva), entonces stas son altamente regulares. De hecho, son monomios del tipo $f(x)=ax$, y tienen, por tanto, el ms alto grado de suavidad posible: son funciones analticas. Surge, por tanto, la siguiente pregunta natural: Àhay otras ecuaciones funcionales que tengan una propiedad similar?. En este artculo vamos a explicar algunos resultados relacionados con esta cuestin para el caso de otra ecuacin funcional clsica: la ecuacin de Frchet.

\section{El teorema de Fr\'{e}chet clsico}

En 1909 el matemtico francs Maurice Frchet \cite{frechet} demostr que, entre todas las funciones continuas $f:\mathbb{R}\to\mathbb{R}$, los polinomios de grado $\leq m$ se pueden caracterizar como las soluciones de una ecuacin funcional que generaliza de forma natural a la ecuacin de Cauchy. Concretamente, demostr el siguiente resultado:

\begin{teo}[Frchet, 1909] \label{Pre_Teo_Frechet}
Consideremos el operador
\begin{equation*}
\begin{array}{l}
\mathcal{F}_{m+1}(f)(x_{1},\cdots ,x_{m+1}) = f(x_{1}+x_{2}+\cdots +x_{m+1})\\
+\sum_{t=1}^{m}(-1)^{t}\sum_{\{i_{1},\cdots ,i_{m+1-t}\}\in
\mathcal{P}_{t}(m+1)}f(x_{i_{1}}+\cdots
+x_{i_{m+1-t}})+(-1)^{m+1}f(0)\text{,}
\end{array}
\end{equation*}
donde $x_1,x_2,\cdots,x_{m+1}$ son variables reales y  
\begin{equation*}
\mathcal{P}_{t}(m+1)=\{A\subset \{1,2,...,m+1\}:\#A=m+1-t\},\text{ }%
t=1,2,...,m.
\end{equation*}
Si  $f:\mathbb{R\rightarrow R}$ es una funcin continua, entonces se tiene que $f$   es un polinomio algebraico de grado menor o igual a $m$ (es decir, $f(x)=a_0+a_1x+\cdots+a_mx^m$ para ciertas constantes $\{a_i\}_{i=0}^m\subseteq \mathbb{R}$ y para todo $x\in\mathbb{R}$) si y solo si 
la funcin $\mathcal{F}_{m+1}(f)$ se anula identicamente en $\mathbb{R}^{m+1}$.
\end{teo}

A continuacin vamos a repetir los argumentos utilizados originalmente por Frchet para su prueba del Teorema  \ref{Pre_Teo_Frechet}. Existen otras demostraciones posibles y, de hecho, en este artculo vamos a presentar tambin varias demostraciones nuevas del mismo resultado. Aunque la importancia de la ecuacin de Frchet (ni de ninguna otra ecuacin funcional que se pueda considerar) no es comparable a la de la ecuacin de Cauchy, lo cierto es que sta ha suscitado, a lo largo del tiempo, el inters de numerosos matemticos, dando lugar a algunos de los artculos que, en nuestra opinin, podran considerarse ms hermosos de la teora de ecuaciones funcionales.

Por comodidad, para evitar un uso excesivo de notacin y frmulas engorrosas, vamos a demostrar, antes de abordar la prueba del Teorema de Frchet, algunos resultados tcnicos que luego sern de utilidad. Comenzamos con la definicin del operador en diferencias progresivas de orden $n$:

\begin{concepto}
Dada una funcin $f:\mathbb{R}\to\mathbb{R}$, definimos los operadores $$\Delta_hf(x)=f(x+h)-f(x) \ \text{ (para } h,x\in\mathbb{R}\text{),}$$ y $$\Delta_{h_1h_2\cdots h_s}f(x)=\Delta_{h_1}\left(\Delta_{h_2\cdots h_s}f\right)(x),\  s=2,3,\cdots,\  (x,h_1,h_2,\cdots,h_s)\in\mathbb{R}^{s+1}.$$ 
Adems, si $h_1=h_2=\cdots=h_s=h$, usamos la notacin $\Delta_h^sf(x)$ para  $\Delta_{h,h,\cdots ,h}f(x)$.
\end{concepto}

\begin{lem}$$\Delta_{h_1h_2}=\Delta_{h_1+h_2}-\Delta_{h_1}-\Delta_{h_2}=\Delta_{h_2h_1}$$
\end{lem}

\noindent \textbf{Demostracin. } Para probar la primera igualdad, basta hacer los clculos:
\begin{eqnarray*}
\Delta_{h_1h_2}f(x) &=& \Delta_{h_1}(f(x+h_2)-f(x)) \\
&=& f(x+h_2+h_1)-f(x+h_1)-f(x+h_2)+f(x)\\
&=& f(x+h_2+h_1)-f(x)+f(x)- f(x+h_1)+f(x)-f(x+h_2)\\
&=& \Delta_{h_1+h_2}f(x)-\Delta_{h_1}f(x)-\Delta_{h_2}f(x).
\end{eqnarray*}
La segunda igualdad es consecuencia de que intercambiar el orden de $h_1,h_2$ en el segundo miembro de la expresin anterior, no cambia nada. {\hfill $\Box$}

\begin{lem} Se tiene que  $$ \Delta_{h_1h_2\cdots h_s}f(x)=\sum_{\varepsilon_1,\varepsilon_2,\cdots,\varepsilon_s=0}^1(-1)^{s-(\varepsilon_1+\varepsilon_2+\cdots+\varepsilon_s)}f(x+\varepsilon_1h_1+\varepsilon_2h_2+\cdots +\varepsilon_sh_s).$$
En particular, $$ \mathcal{F}_{m+1}(f)(x_{1},\cdots ,x_{m+1}) = \Delta_{x_{1}\cdots x_{m+1}}f(0).$$
\end{lem}

\noindent \textbf{Demostracin. } Este resultado se demuestra por un proceso rutinario de induccin. En efecto, para $s=2$ pasos el resultado es trivial. Supongamos que es cierto para $s$ pasos y veamos qu sucede cuando tomamos $s+1$ pasos:
\begin{eqnarray*}
&\ &  \Delta_{h_{s+1}}(\Delta_{h_1h_2\cdots h_s}f(x)) \\
&=& \sum_{\varepsilon_1,\varepsilon_2,\cdots,\varepsilon_s=0}^1(-1)^{s-(\varepsilon_1+\varepsilon_2+\cdots+\varepsilon_s)}f(x+\varepsilon_1h_1+\varepsilon_2h_2+
\cdots+\varepsilon_sh_s+h_{s+1}) \\
&\ & \ \ -\sum_{\varepsilon_1,\varepsilon_2,\cdots,\varepsilon_s=0}^1(-1)^{s-(\varepsilon_1+\varepsilon_2+\cdots+\varepsilon_s)}f(x+\varepsilon_1h_1+\varepsilon_2h_2+\cdots+\varepsilon_sh_s) \\
&=& \sum_{\varepsilon_1,\varepsilon_2,\cdots,\varepsilon_s=0}^1(-1)^{s+1-(\varepsilon_1+\varepsilon_2+\cdots+\varepsilon_s+1)}f(x+\varepsilon_1h_1+\varepsilon_2h_2+\cdots+\varepsilon_sh_s+h_{s+1}) \\
&\ & \ \ +\sum_{\varepsilon_1,\varepsilon_2,\cdots,\varepsilon_s=0}^1(-1)^{s+1-(\varepsilon_1+\varepsilon_2+\cdots+\varepsilon_s)}f(x+\varepsilon_1h_1+\varepsilon_2h_2+\cdots+\varepsilon_sh_s) \\
&=& \sum_{\varepsilon_1,\varepsilon_2,\cdots,\varepsilon_{s+1}=0}^1(-1)^{s+1-(\varepsilon_1+\varepsilon_2+\cdots+\varepsilon_{s+1})}f(x+\varepsilon_1h_1+\varepsilon_2h_2+\cdots+\varepsilon_{s+1}h_{s+1}),
\end{eqnarray*}
que es lo que buscbamos.  La segunda afirmacin del lema es consecuencia inmediata de la frmula que acabamos de probar y de la definicin del operador $ \mathcal{F}_{m+1}$. {\hfill $\Box$}

\noindent\textbf{Demostracin del Teorema de Frchet.}  Hacemos la prueba por induccin sobre $m$.  Si $m=1$, entonces
\[
\mathcal{F}_{2}(f)=f(x_1+x_2)-f(x_1)-f(x_2)+f(0),
\]
por lo que, si $\mathcal{F}_{2}(f)$ se anula idnticamente, entonces, tomando $g(x)=f(x)-f(0)$, se comprueba que
\begin{eqnarray*}
g(x_1+x_2)-g(x_1)-g(x_2)&=& f(x_1+x_2)-f(0)-f(x_1)+f(0)-f(x_2)+f(0)\\
&=& f(x_1+x_2)-f(x_1)-f(x_2)+f(0)=0,
\end{eqnarray*}
por lo que $g$ es una solucin continua de la ecuacin de Cauchy, de modo que $g(x)=a_1x$ para cierto nmero real $a_1$, y, en consecuencia, $f(x)=a_0+a_1x$, con $a_0=f(0)$. Esto demuestra el resultado para $m=1$.

Supongamos que el teorema es cierto para $m-1$ y consideremos una solucin continua $f$ de la ecuacin  
$\mathcal{F}_{m+1}(f)=0$. Introducimos la funcin auxiliar 
\begin{equation*}
\varphi (x)=f(x+x_{m+1})-f(x)-f(x_{m+1})+f(0).
\end{equation*}
No es difcil comprobar que si  $\mathcal{F}_{m+1}(f)=0$ entonces, para cada constante fija $x_{m+1}$, se tiene que 
$\mathcal{F}_{m}(\varphi)=0$. En efecto, es evidente que 
\[
\varphi(x)=\Delta_{x_{m+1}}f(x)-\Delta_{x_{m+1}}f(0)=\Delta_{x_{m+1}}(\Delta_xf(0))=\Delta_{x_{m+1}x}f(0),
\]
por lo que
\begin{eqnarray*}
\mathcal{F}_{m}(\varphi)(x_{1},\cdots ,x_{m+1}) &=& \Delta_{x_1x_2\cdots x_m}\Delta_{x_{m+1}}(\Delta_xf(0))\\
&=& 
 \Delta_{x_1x_2\cdots x_{m+1}}(\Delta_xf(0))\\
&=& \Delta_x( \Delta_{x_1x_2\cdots x_{m+1}}f)(0)=0
\end{eqnarray*}
Se sigue que podemos utilizar la hiptesis de induccin para la funcin $\varphi$. Es decir, ahora sabemos que 
 $\varphi (x)\in \Pi_{m-1} $. 

Por otra parte, la funcin 
\begin{equation*}
Q(x,y)=f(x+y)-f(x)-f(y)+f(0)
\end{equation*}
es evidentemente simtrica (es decir: $Q(x,y)=Q(y,x)$ para todo par $(x,y)$), lo que nos conduce a concluir que  $Q(x,y)$ 
es un polinomio en las variables $x,y$ y su grado es menor o igual a $m-1$ en cada una de estas variables.

Dividimos el resto de la demostracin en varias etapas:

\begin{itemize}
\item  Primero observamos que $Q(x,y)$ satisface la ecuacin funcional
\begin{equation*}
Q(x,y)+Q(x+y,z)=f(x+y+z)-f(x)-f(y)-f(z)+2f(0),
\end{equation*}
por lo que  $Q(x,y)+Q(x+y,z)$ es una funcin simtrica con respecto a las variables 
 $x, y, z$.

\item  Descomponemos $Q(x,y)$ como suma de sus componentes homogneas,  
\begin{equation*}
Q=Q_{0}+Q_{1}+\cdots +Q_{r}\text{,}
\end{equation*}
donde $Q_{i}$ es un polinomio homogneo de grado  $i$, y $i=0,\cdots ,r$.
Como esta descomposicin es nica para cualquier polinomio en un nmero finito de variables, y  como, para cada $i$, el polinomio $Q_{i}(x,y)+Q_{i}(x+y,z)$ es homogneo de grado $i$,  concluimos que los polinomios 
$Q_{i}(x,y)+Q_{i}(x+y,z)$ son funciones simtricas respecto de las variables $x,y,z$, por serlo el polinomio $Q(x,y)+Q(x+y,z)$.

\item  Como $Q_{i}(x,y)$ es homogneo de grado $i$, sabemos que admite una expresin del tipo:
\begin{equation*}
Q_{i}(x,y)=a_{0}x^{i}+a_{1}x^{i-1}y+\cdots +a_{i-1}xy^{i-1}+a_{i}y^{i},
\end{equation*}
por lo que, si utilizamos la simetra de la funcin  $Q_{i}(x,y)+Q_{i}(x+y,z)$ conjuntamente con la frmula del Binomio de Newton, obtenemos un conjunto de identidades sobre los coeficientes 
$\{a_{t}\}_{t=0}^{i}$ que, si se usan apropiadamente, nos garantizan la existencia de cierta constante
 $A_{i}$ tal que 
\begin{equation*}
Q_{i}(x,y)=A_{i}\left( (x+y)^{i}-x^{i}-y^{i}\right) ,
\end{equation*}
para $i=1,2,...,r.$ (Dejamos para el lector la comprobacin, usando un proceso de induccin, de esta afirmacin).

En particular, los monomios  $x^{i}$, $y^{i}$ no aparecen en la expresin del polinomio 
 $Q_{i}(x,y)$, por lo que, si $Q=Q_{0}+\cdots +Q_{r}$ tiene grado  
 $\leq m-1$ en cada una de las variables  $x$ e $y$, entonces podemos afirmar que $r\leq
m$ y
\begin{equation*}
Q(x,y)=\sum_{i=2}^{m}A_{i}\left( (x+y)^{i}-x^{i}-y^{i}\right)
=R(x+y)-R(x)-R(y),
\end{equation*}
donde $R(x)=\sum_{i=2}^{m}A_{i}x^{i}$.

\item  Consideramos ahora la funcin  $S(x)=f(x)-R(x)-f(0)$. Esta funcin satisface la ecuacin
\begin{equation*}
S(x+y)-S(x)-S(y)=0 \text{ para todo }x,y\in \mathbb{R}\text{,}
\end{equation*}
por lo que existe una cierta constante $a\in\mathbb{R}$ tal que $S(x)=ax$. Se sigue que 
$f(x)=f(0)+ax+R(x)\in \Pi_{m}$, lo que finaliza la demostracin. 
\end{itemize}

{\hfill $\Box$}

Es importante observar que la demostracin original de Frchet se puede repetir sin grandes cambios bajo la hiptesis, mucho ms dbil, de que la funcin $f$ es acotada en un conjunto $A\subseteq \mathbb{R}$ con medida de Lebesgue positiva. Para ello, basta utilizar un teorema debido a Kurepa $\cite{kurepa}$ que generaliza el Teorema de Steinhaus $\cite{steinhaus}$. Concretamente, dicho resultado garantiza que la funcin $\lambda(x)=|A\cap (A-\{x\})|$ es continua, por lo que, si $|A|=\lambda(0)>0$, entonces $\lambda(x)>0$ para todo $|x|<\varepsilon$, para cierto $\varepsilon>0$. Esto implica que las funciones del tipo $g(x)=f(x+a)-f(x)$ mantienen la propiedad de estar acotadas en un conjunto de medida positiva, siempre que $|a|$ sea suficientemente pequeo. 
Supongamos ahora que  $f$ es una solucin de  $\mathcal{F}_{m+1}(f)=0$ y $\sup_{x\in A}|f(x)|\leq M<\infty$ para cierto conjunto $A$ con $|A|>0$. Entonces la funcin  
\[
\varphi (x)=f(x+x_{m+1})-f(x)-f(x_{m+1})+f(0)
\]
tambin est acotada en el conjunto $A$  (siempre que $|x_{m+1}|$ sea suficientemente pequeo) y, como se prob anteriormente, satisface $\mathcal{F}_{m}(\varphi)=0$, por lo que podemos aplicarle a ella la hiptesis de induccin. Esto nos conduce a que $\varphi\in \Pi_{m-1}$. A continuacin, es evidente que los argumentos dados relativos a la funcin $Q(x,y)=f(x+y)-f(x)-f(y)+f(0)$ no requieren cambios. Finalmente, como $R(x)$ es un polinomio algebraico ordinario, la funcin $S(x)=f(x)-R(x)-f(0)$ est acotada en $A$ y, al aplicar a ella el Teorema de Kormes ($\text{Teorema } \ref{kormes}$), tenemos que $S(x)=ax$ para cierta constante $a$, que es lo que buscbamos.

La ecuaci\'{o}n $\Delta_{h_1h_2\cdots h_{m+1}}f(x)=0 $ se puede estudiar para funciones $f:X\to Y$  cuando $X, Y$ son un par de espacios vectoriales sobre  $\mathbb{Q}$, y las variables  $x,h_1,\cdots,h_{m+1}$ representan elementos de $X$, 
\begin{equation}\label{Pre_fregeneral}
\Delta_{h_1h_2\cdots h_{m+1}}f(x)=0 \ \ (x,h_1,h_2,\dots,h_{m+1}\in X).
\end{equation}
En este contexto, las soluciones de \eqref{Pre_fregeneral} son  las funciones de la forma $f(x)=A_0+A_1(x)+\cdots+A_m(x)$, donde  $A_0$ es una constante y 
$A_k(x)=A^k(x,x,\cdots,x)$ para cierta funci\'{o}n  $k$-aditiva sim\'{e}trica  $A^k:X^k\to Y$ (decimos que   $A_k$ es la diagonalizaci\'{o}n de  $A^k$) y para $k=0,1,\cdots,m$. En particular, si $x\in X$ y 
$r\in\mathbb{Q}$, entonces $f(rx)=A_0+rA_1(x)+\cdots+r^mA_m(x)$. Adem\'{a}s, se sabe que $f:X\to Y$ satisface  \eqref{Pre_fregeneral} si y solo si es una soluci\'{o}n de la ecuaci\'{o}n funcional 
\begin{equation}\label{Pre_frepasofijo}
\Delta_{h}^{m+1}f(x):=\sum_{k=0}^{m+1}\binom{m+1}{k}(-1)^{s-k}f(x+kh)=0 \ \ (x,h\in X).
\end{equation}
Una demostraci\'{o}n de este resultado se sigue de un conocido teorema, debido a  Djokovi\'{c} \cite{Dj} (ver tambi\'{e}n  \cite[Theorem 7.5, p\'{a}g. 160]{HIR}, \cite[Theorem 15.1.2., p\'{a}g. 418]{kuczma}), el cual establece que los operadores $\Delta_{h_1 h_2\cdots h_s}$ verifican la ecuaci\'{o}n 
\begin{equation}\label{Pre_igualdad}
\Delta_{h_1\cdots h_s}f(x)=
\sum_{\epsilon_1,\dots,\epsilon_s=0}^1(-1)^{\epsilon_1+\cdots+\epsilon_s}
\Delta_{\alpha_{(\epsilon_1,\dots,\epsilon_s)}(h_1,\cdots,h_s)}^sf(x+\beta_{(\epsilon_1,\dots,\epsilon_s)}(h_1,\cdots,h_s)),
\end{equation}
donde $$\alpha_{(\epsilon_1,\dots,\epsilon_s)}(h_1,\cdots,h_s)=(-1)\sum_{r=1}^s\frac{\epsilon_rh_r}{r}$$ y 
$$\beta_{(\epsilon_1,\dots,\epsilon_s)}(h_1,\cdots,h_s)=\sum_{r=1}^s\epsilon_rh_r.$$  
Si tenemos ahora en cuenta la equivalencia de las ecuaciones de Frchet con paso fijo y con paso variable, y utilizamos los argumentos expuestos tras la prueba del teorema de Frchet, habremos demostrado el siguiente resultado:
\begin{teo}\label{frechet_kormes} Supongamos que $f:\mathbb{R}\to\mathbb{R}$ es solucin de la ecuacin de Frchet $\Delta_h^{m+1}f(x)=0$ y $f$ es acotada en $A\subseteq \mathbb{R}$ para cierto conjunto de medida positiva, $|A|>0$. Entonces $f\in \Pi_m$.
\end{teo}
El Teorema \ref{frechet_kormes} ha sido demostrado varias veces, con mtodos muy diversos. De hecho, el argumento que hemos presentado aqu para su demostracin, aunque est basado en la prueba original de Frchet, no se ha publicado hasta la fecha -al menos, por lo que alcanza a nuestro conocimiento del tema- y puede considerarse, en consecuencia, nuevo. Sin embargo, nuestra prueba depende tambin de la herramienta tpica con la que se ha demostrado este resultado anteriormente: el Teorema de Kurepa. Por ejemplo, una demostracin de este tipo la podemos encontrar el la monografa de  Sz\'{e}kelyhidi \cite{laszlo1}. Otras pruebas se pueden encontrar, por ejemplo, en  \cite{ciesielski}, \cite{ger}, \cite{mckiernan}. En particular, Ciesielski  \cite{ciesielski} demostr que, si $f:\mathbb{R}\to\mathbb{R}$ es solucin de la siguiente desigualdad funcional
\begin{equation*}
\Delta _{h}^{n+1}f(x)\geq 0\text{, para todo }x,h\in \mathbb{R}
\end{equation*}
y est acotada en un conjunto $A $ con medida de Lebesgue positiva, entonces $f$ es una funcin continua en toda la recta real. Se sigue que las soluciones de la ecuacin de Frchet que no son polinomios ordinarios no pueden estar acotadas en ninguno de estos conjuntos $A$. 

\section{El teorema de Montel clsico}

En un artculo de 1937, Paul Montel \cite{montel} demostr un resultado relacionado con la ecuacin funcional de Frchet que, muy probablemente, result inesperado. En vez de centrarse en las propiedades mnimas de regularidad necesarias para que la funcin $f$, solucin de la ecuacin en cuestin, se viera forzada a ser un polinomio algebraico ordinario, se ocup de conocer cuntos pasos $h$ son estrictamente necesarios para garantizar que una funcin continua que satisface $\Delta_{h_k}^{m+1}f(x)=0$ (para todo $x$ y para unos pocos valores de $k$), es por necesidad un polinomio. Como l mismo explicaba en su artculo, la motivacin principal para este estudio no proceda de los resultados publicados por Frchet en 1909, sino de un teorema muy anterior, debido a C.G. Jacobi \cite{Jacobi} y publicado en 1834, sobre los periodos de las funciones meromorfas. En efecto, si tomamos $n=1$, decir que $\Delta_h^1f(x)=0$ para todo $x$ es lo mismo que decir que $f$ es una funcin peridica y $h$ es uno de sus periodos. Jacobi haba demostrado, mediante una construccin explcita, que existen funciones $f:\mathbb{C}\to\widehat{\mathbb{C}}$ que son meromorfas, no son constantes, y admiten dos periodos independientes.  Estas funciones reciben el nombre de funciones doblemente peri\'{o}dicas (o el\'{\i}pticas) y son especialmente importantes para la variable compleja  \cite{JS}. No existen funciones holomorfas no constantes que sean doblemente peridicas, porque stas deben ser funciones acotadas y, en tal caso, el Teorema de Liouville garantiza que son funciones constantes.  Adems, si una funcin meromorfa admite tres periodos independientes (ya diremos lo que significa esto), entonces forzosamente esta funcin es constante. Montel se preguntaba si algo similar podra suceder con los ``periodos generalizados'' $h$ que aparecen en la ecuacin de Frchet  $\Delta_{h}^{m+1}f=0$.    
En efecto, debemos observar que cuando $m=0$, los periodos $h$ forman un subgrupo aditivo de $\mathbb{R}$, cosa que, como veremos posteriormente, no sucede necesariamente para el resto de valores de $m$ y que tiene el efecto de hacer el problema muy sencillo cuando $m=0$, pero complicado cuando $m\geq 1$.

Comenzamos, pues, realizando algunas observaciones sencillas, formuladas en un contexto ms general que el originalmente planteado por Jacobi, sobre el caso $m=0$. 
Dados un grupo conmutativo  $(G,+)$, un conjunto no vac\'{\i}o $Y$, y una funci\'{o}n $f:G\to Y$, consideramos el conjunto de los periodos de $f$, $\mathfrak{P}_0(f)=\{g\in G:f(w+g)=f(w)\text{ 
para todo  } w\in G\}$. Obviamente,   $\mathfrak{P}_0(f)$ es siempre un subgrupo de   $G$.  En algunos casos especiales, estos grupos son conocidos y, de hecho, poseen una hermosa estructura. Por ejemplo, el Teorema de Jacobi de 1834 establece que  si $f:\mathbb{C}\to\widehat{\mathbb{C}}$ es una funci\'{o}n meromorfa no constante, definida sobre los n\'{u}meros complejos, entonces    $\mathfrak{P}_0(f)$ es un subgrupo discreto de   $(\mathbb{C},+)$, lo cual reduce el estudio de estos conjuntos a los siguientes tres casos:  
$\mathfrak{P}_0(f)=\{0\}$,  $\mathfrak{P}_0(f)=\{nw_1:n\in\mathbb{Z}\}$  para cierto n\'{u}mero complejo $w_1\neq 0$, o $\mathfrak{P}_0(f)=\{n_1w_1+n_2w_2:(n_1,n_2)\in \mathbb{Z}^2\}$ 
para un par de n\'{u}meros complejos $w_1,w_2$ que satisfacen $w_1w_2\neq 0$ y $w_1/w_2\not\in\mathbb{R}$. En particular, estas funciones no pueden tener tres periodos independientes y, adem\'{a}s, existen ejemplos de funciones meromorfas $f:\mathbb{C}\to\widehat{\mathbb{C}}$  con dos periodos independientes  $w_1,w_2$ siempre que $w_1/w_2\not\in\mathbb{R}$. De forma an\'{a}loga, si la funci\'{o}n $f:\mathbb{R}\to\mathbb{R}$ es continua y no constante, entonces no admite dos periodos que sean linealmente independientes sobre $\mathbb{Q}$.  Estos resultados pueden formularse en t\'{e}rminos de ecuaciones funcionales, puesto que  $h$ es un periodo de $f:G\to Y$ si y solo si $f$ resuelve la ecuaci\'{o}n funcional 
$\Delta_hf(x)=0 \ \ (x\in G).$ Por tanto, el Teorema de Jacobi se puede reformular como un teorema sobre ecuaciones funcionales en el que se establece que las funciones constantes son precisamente  las funciones meromorfas  $f:\mathbb{C}\to\widehat{\mathbb{C}}$  que resuelven un sistema de ecuaciones  del tipo 
\begin{equation}\label{Pre_JC}
\Delta_{h_1}f(z)=\Delta_{h_2}f(z)=\Delta_{h_3}f(z)=0 \ \ (z\in \mathbb{C})
\end{equation}
para tres periodos independientes $\{h_1,h_2,h_3\}$ (i.e., $h_1\mathbb{Z}+ h_2\mathbb{Z}+h_3\mathbb{Z}$ es un subconjunto denso de $\mathbb{R}$). En el caso real, el teorema afirma que, si $h_1h_2\neq0$ y $h_1/h_2\not\in\mathbb{Q}$, entonces la funci\'{o}n continua   $f:\mathbb{R}\to\mathbb{R}$ es constante si y solo si es soluci\'{o}n del sistema de ecuaciones funcionales
\begin{equation} \label{Pre_JR}
\Delta_{h_1}f(x)=\Delta_{h_2}f(x)= 0\ (x\in \mathbb{R}).
\end{equation}
Paul Montel  sustituy\'{o}
en las ecuaciones  $\eqref{Pre_JC},\eqref{Pre_JR}$ el operador en diferencias finitas de orden uno, $\Delta_h$,  por el operador en diferencias progresivas de orden superior,  $\Delta^{m+1}_h$ y prob\'{o} que las ecuaciones resultantes son apropiadas para la caracterizaci\'{o}n de los polinomios algebraicos. Concretamente, demostr\'{o} el siguiente resultado:

\begin{teo}[Montel]\label{Pre_montel_1v} Las siguientes afirmaciones son ciertas:
\begin{itemize}
\item[$(i)$]  Si $\{h_1,h_2\}\subset \mathbb{R}$ son tales que $h_1h_2\neq 0$ y  $h_1/h_2\not\in\mathbb{Q}$, entonces la funci\'{o}n continua  $f:\mathbb{R}\to\mathbb{R}$ es un polinomio algebraico con coeficientes reales y grado  $\leq m$ (i.e., $f(x)=a_0+a_1x+\cdots+a_mx^m$ con los coeficientes $a_i$ nmeros reales para todo $i$) si y solo si resuelve el sistema de ecuaciones funcionales 
\begin{equation} \label{Pre_JRM}
\Delta_{h_1}^{m+1}f(x)=\Delta_{h_2}^{m+1}f(x)= 0\ (x\in \mathbb{R}).
\end{equation}
\item[$(ii)$] Supongamos que $f:\mathbb{C}\to\mathbb{C}$ es una funci\'{o}n holomorfa que resuelve un sistema de ecuaciones funcionales de la forma 
\begin{equation}\label{Pre_JCM}
\Delta_{h_1}^{m+1}f(z)=\Delta_{h_2}^{m+1}f(z)=\Delta_{h_3}^{m+1}f(z)=0 \ \ (z\in \mathbb{C})
\end{equation}
para tres periodos independientes  $\{h_1,h_2,h_3\}$. Entonces $f(z)=a_0+a_1z+\cdots+a_mz^m$ es un polinomio algebraico con coeficientes complejos y grado   $\leq m$. 
\end{itemize}
\end{teo}

Veamos cmo demostr Montel el teorema anterior. Para ello, recurriremos a varios resultados tcnicos intermedios.

\begin{lem} \label{Pre_lemM1} Supongamos que $h_1,h_2\in\mathbb{R}$ son tales que $h_1/h_2\not\in\mathbb{Q}$. Si $f:\mathbb{R}\to\mathbb{R}$ es continua y 
$\Delta_{h_1}f(x)=\Delta_{h_2}f(x)=0$ para todo $x\in\mathbb{R}$, entonces $f(x)=c$ es una funcin constante.
\end{lem}

\noindent \textbf{Demostracin. } El resultado es consecuencia inmediata de la continuidad de $f$ y de que si  $h_1/h_2\not\in\mathbb{Q}$, entonces $h_1\mathbb{Z}+h_2\mathbb{Z}$ es un subconjunto denso de $\mathbb{R}$ y todos sus elementos son periodos de $f$. {\hfill $\Box$}

\begin{lem} \label{Pre_lemM2} Supongamos que $h_1,h_2\in\mathbb{R}$ son tales que $h_1/h_2\not\in\mathbb{Q}$. Si $f:\mathbb{R}\to\mathbb{R}$ es continua y 
\begin{equation} \label{Pre_lineal}
\Delta_{h_1}f(x)=c_1,\  \Delta_{h_2}f(x)=c_2 \text{ para todo } x\in\mathbb{R} \text{ y ciertas constantes }c_1,c_2,
\end{equation} 
entonces $f(x)=ax+b$ para ciertas constantes $a,b$. Adems, si $c_1c_2=0$, entonces $f(x)=c$ es una funcin constante.
\end{lem}

\noindent \textbf{Demostracin. } Supongamos que $f$ satisface \eqref{Pre_lineal}, y consideremos, para cada $h$, la funcin $g_h(x)=\Delta_hf(x)=f(x+h)-f(x)$. Entonces
$$\Delta_{h_i}g_h(x)=\Delta_{h_i}\Delta_hf(x)=\Delta_h\Delta_{h_i}f(x)=0,\ i=1,2,$$
por lo que, aplicando el Lema \ref{Pre_lemM1}, concluimos que $g_h(x)=g(0)=f(h)-f(0)$, para todo $x,h\in\mathbb{R}$. En otras palabras, la funcin $f$ es continua y satisface la ecuacin funcional
\[
f(x+h)=f(x)+f(h)-f(0),\ \text{ para todo } x,h\in\mathbb{R}.
\]
Esto nos conduce directamente a que $f(x)=ax+b$ para todo $x\in\mathbb{R}$ y ciertas constantes $a,b$. 

Para demostrar la segunda parte del lema basta observar que si $f(x)=ax+b$ entonces $\Delta_hf(x)=ah$ para todo $x,h\in\mathbb{R}$. Ahora bien, si $c_1c_2=0$, entonces $c_1=0$ o $c_2=0$ y en ambos casos lo que estamos diciendo es que $a=0$ y, por tanto, $f(x)=b$ es una funcin constante. 
{\hfill $\Box$}

\begin{lem} \label{Pre_lemM3} Supongamos que $\{h_1,h_2,h_3\}\subseteq \mathbb{C}$ son tres periodos independientes 
para la funcin continua $f:\mathbb{C}\to\mathbb{C}$. Entonces $f(z)=c$ para cierta constante $c\in\mathbb{C}$.
\end{lem}
\noindent \textbf{Demostracin. } El resultado es consecuencia inmediata de la continuidad de $f$ y de que si los periodos $\{h_1,h_2,h_3\}\subseteq \mathbb{C}$ son independientes, entonces $h_1\mathbb{Z}+h_2\mathbb{Z}+h_2\mathbb{Z}$ es un subconjunto denso de $\mathbb{C}$.
 
{\hfill $\Box$}

\begin{lem}\label{Pre_lemM4} Supongamos que $\{h_1,h_2,h_3\}\subseteq \mathbb{C}$ son independientes, y que  $f:\mathbb{C}\to\mathbb{C}$ es holomorfa y 
$\Delta_{h_k}f(z)=c_k$, $k=1,2,3$, para ciertas constantes $c_1,c_2,c_3\in\mathbb{C}$. Entonces $f(z)=az+b$ para ciertas constantes $a,b\in\mathbb{C}$. Adems, si $c_1c_2c_3=0$, entonces $f(z)=b$ es una funcin constante. 
\end{lem}

\noindent \textbf{Demostracin. } En efecto, como estamos asumiendo que $f(z)$ puede derivarse, y que la derivada es una funcin continua, se tiene que, para cada $k\in\{1,2,3\}$,
\begin{eqnarray*}
\Delta_{h_k}f'(z) &=&  \Delta_{h_k}\lim_{h\to 0}\frac{f(z+h)-f(z)}{h}\\
&=&  \lim_{h\to 0}\frac{\Delta_{h_k}f(z+h)-\Delta_{h_k}f(z)}{h}\\
&=& \lim_{h\to 0} \frac{0}{h}=0.
\end{eqnarray*}
Por tanto, el Lema \ref{Pre_lemM1} implica que $f'(z)$ es una funcin constante. Es decir, $f(z)=az+b$ para ciertos valores $a,b\in\mathbb{C}$. 

La segunda parte del lema es consecuencia de que si $f(x)=ax+b$ entonces $\Delta_hf(z)=ah$ para todo $z,h\in\mathbb{C}$ pues, si $c_1c_2c_3=0$, entonces $c_1=0$ o $c_2=0$ o $c_3=0$ y  en todos estos casos se concluye que  $a=0$ y, por tanto, $f(x)=b$ es una funcin constante. 

{\hfill $\Box$}

\noindent \textbf{Demostracin del Teorema de Montel.} Aunque las ideas subyacentes son muy similares, hacemos una prueba separada para cada apartado del teorema.

\noindent $(i)$. Supongamos que 
$h_1,h_2\in\mathbb{R}$ son tales que $h_1h_2\neq 0$ y $h_1/h_2\not\in\mathbb{Q}$ y  $f:\mathbb{R}\to\mathbb{R}$ es continua y 
\begin{equation} \label{Pre_polinom}
\Delta_{h_1}^{m+1}f(x) = \Delta_{h_2}^{m+1}f(x)=0 \text{ para todo } x\in\mathbb{R}.
\end{equation} 
Vamos a proceder por induccin sobre $m$. Para $m=0$ el resultado coincide con el Lema \ref{Pre_lemM1}, que ya se demostr. 

Consideremos la funcin 
\[
\varphi_0(x)=\Delta_{h_1}^m\Delta_{h_2}^mf(x).
\]
Entonces $\varphi_0$ es continua y $\Delta_{h_i}\varphi(x)=0$, $i=1,2$. Por tanto, el Lema \ref{Pre_lemM1} garantiza que $\varphi_0(x)=\alpha_0$ para cierta constante $\alpha_0$. Una vez tenemos esta informacin, consideramos la nueva funcin 
\[
\varphi_1(x)=\Delta_{h_1}^m\Delta_{h_2}^{m-1}f(x).
\] 
Entonces $\varphi_1(x)$ es  continua y 
\begin{eqnarray*}
\Delta_{h_1}\varphi_1(x)=\Delta_{h_1}^{m+1}\Delta_{h_2}^{m-1}f(x)= \Delta_{h_2}^{m-1}\Delta_{h_1}^{m+1}f(x)=0\\
\Delta_{h_2}\varphi_1(x)=\Delta_{h_1}^{m}\Delta_{h_2}^{m}f(x)= \varphi_0(x)=\alpha_0.
\end{eqnarray*}
En otras palabras, $\varphi_1(x)$ satisface la ecuacin \eqref{Pre_lineal} con $c_1=0$, $c_2=\alpha_0$, por lo que, aplicando el Lema \ref{Pre_lemM2}, tenemos que 
$\alpha_0=0$ y  $\varphi_1(x)=\alpha_1$  es una funcin constante. 

Este mismo argumento lo podemos repetir inductivamente con las funciones
\[
\varphi_i(x)=\Delta_{h_1}^m\Delta_{h_2}^{m-i}f(x), \text{  } i=0,1,2,\cdots,m.
\] 
En el paso $i$-simo, tendremos que $\varphi_i(x)=\alpha_i$ para cierta constante $\alpha_i$. Entonces, aplicando los operadores $\Delta_{h_k}$ ($k=1,2$) a la funcin 
$\varphi_{i+1}(x)$, tendremos que
\begin{eqnarray*}
\Delta_{h_1}\varphi_{i+1}(x)=\Delta_{h_1}^{m+1}\Delta_{h_2}^{m-(i+1)}f(x)= \Delta_{h_2}^{m-(i+1)}\Delta_{h_1}^{m+1}f(x)=0\\
\Delta_{h_2}\varphi_{i+1}(x)=\Delta_{h_1}^{m}\Delta_{h_2}^{m-i}f(x)= \varphi_i(x)=\alpha_i,
\end{eqnarray*}
por lo que aplicando el Lema \ref{Pre_lemM2}, tendremos que 
$\alpha_i=0$ y  $\varphi_{i+1}(x)=\alpha_{i+1}$  es una funcin constante. Evidentemente, el proceso termina cuando $i=m$, lo que nos lleva a la conclusin de que existe una constante $C_1$ (de hecho, $C_1=\alpha_m$) tal que 
\[
\Delta_{h_1}^mf(x)=C_1,\text{ para todo } x\in\mathbb{R}.
\] 
El mismo argumento, pero intercambiando los papeles de $h_1$ y $h_2$, nos sirve para demostrar que 
\[
\Delta_{h_2}^mf(x)=C_2,\text{ para todo } x\in\mathbb{R},
\] 
para cierta constante $C_2$. 

Consideramos ahora la funcin $g(x)=f(x)-\frac{C_1}{m!h_1^m}x^m$. Es claro que 
\begin{eqnarray*}
\Delta_{h_1}^mg(x) &=&  \Delta_{h_1}^mf(x)- \frac{C_1}{m!h_1^m}\Delta_{h_1}^mx^m=C_1-C_1=0, \\
\Delta_{h_2}^mg(x) &=&  \Delta_{h_2}^mf(x)- \frac{C_1}{m!h_1^m}\Delta_{h_2}^mx^m=C_2-C_1(\frac{h_2}{h_1})^m=C.
\end{eqnarray*}
Por tanto, si definimos la funcin 
\[
\psi(x)=\Delta_{h_1}^{m-1}\Delta_{h_2}^{m-1}g(x),
\]
entonces 
\begin{eqnarray*}
\Delta_{h_1}\psi(x) &=& \Delta_{h_1}^{m}\Delta_{h_2}^{m-1}g(x)= \Delta_{h_2}^{m-1}\Delta_{h_1}^{m}g(x)=0\\
\Delta_{h_2}\psi(x) &=& \Delta_{h_1}^{m-1}\Delta_{h_2}^{m}g(x)= 0,
\end{eqnarray*}
y el Lema \ref{Pre_lemM1} garantiza que $\psi(x)=\beta_0$ es una funcin constante. Los mismos argumentos que utilizamos antes con la funcin $f$, cambiando el valor $m+1$ por el valor $m$, nos sirven ahora para probar que $g$ satisface las relaciones
\[
\Delta_{h_1}^{m-1}g(x)=C_1^* \text{ y } \Delta_{h_2}^{m-1}g(x)=C_2^*, \text{ para todo } x\in\mathbb{R} \text{ y ciertas constantes } C_1^*,C_2^*.
\] 
Se sigue que 
\[
\Delta_{h_1}^{m}g(x)= \Delta_{h_2}^{m}g(x)=0, \text{ para todo } x\in\mathbb{R}, 
\]
y, por tanto,  la hiptesis de induccin nos garantiza que $g(x)=a_0+a_1x+\cdots+a_{m-1}x^{m-1}$ para todo $x\in\mathbb{R}$ y ciertos nmeros reales $a_0,\cdots,a_{m-1}$. Esto concluye la demostracin de $(i)$, pues $f(x)=g(x)+ a_mx^m$, con $a_m=\frac{C_1}{m!h_1^m}$.

\noindent $(ii)$. Este apartado admite una demostracin similar a la utilizada para $(i)$.  La diferencia esencial est en que usamos los Lemas \ref{Pre_lemM3} y \ref{Pre_lemM4}, y que la hiptesis de holomorfa simplifica el argumento. La complicacin proviene de que necesitamos usar tres periodos en vez de dos. Razonamos, como es natural, por induccin sobre $m$. El caso $m=0$ est probado con el Lema \ref{Pre_lemM3}. Suponemos que el resultado es cierto para $m-1$ y, a continuacin, consideramos una funcin holomorfa $f:\mathbb{C}\to\mathbb{C}$ que satisface \eqref{Pre_JCM}. Tomamos $\varphi_0(z)=\Delta_{h_1}^m\Delta_{h_2}^m\Delta_{h_3}^mf(z)$. Entonces $\Delta_{h_k}\varphi_0(z)=0$ para $k=1,2,3$, por lo que, aplicando el Lema \ref{Pre_lemM3}, concluimos que $\varphi_0(z)=\gamma_0$ es una funcin constante. Tomamos ahora, como se hizo en la primera parte de la prueba de (i), 
$\varphi_i(z)=\Delta_{h_1}^m\Delta_{h_2}^m\Delta_{h_3}^{m-i}f(z)$, $i=0,1,\dots,m$. Consideremos $\varphi_1(z)$. Es evidente que 
\[
\Delta_{h_1}\varphi_1(z)=\varphi_0(z)=\gamma_0,\text{ y } \Delta_{h_2}\varphi_1(z)=\Delta_{h_3}\varphi_1(z)=0,
\]
por lo que, si usamos el Lema \ref{Pre_lemM4}, concluimos que $\varphi_1(z)=\gamma_1$ es una funcin constante y $\gamma_0=0$. El mismo argumento utilizado en la primera parte de la prueba de $(i)$, nos dice que las funciones $\varphi_i$ satisfacen  $\varphi_i(z)= \Delta_{h_1}^m\Delta_{h_2}^m\Delta_{h_3}^{m-i}f(z)=0$ para $i=0,1,\cdots,m-1$ y $\varphi_m(z)=\Delta_{h_1}^m\Delta_{h_2}^mf(z)=C$ para cierta constante $C$. 
Este tipo de argumento lo podemos ahora aplicar a la funcin $\psi(z)=  \Delta_{h_1}^m\Delta_{h_2}^{m-1}\Delta_{h_3}^{m-1}f(z)$, y concluir que, en realidad, todas las funciones 
\[
 \Delta_{h_1}^m\Delta_{h_2}^{m-1}\Delta_{h_3}^{m-i}f(z), \ i=1,2,\cdots,m-2
\]
se anulan idnticamente, y la funcin $ \Delta_{h_1}^m\Delta_{h_2}^{m-1}f(z)$ es una constante. Nuevamente reiteramos el argumento -tantas veces como sea necesario- y obtendremos que $\Delta_{h_1}^mf(z)=c_1$ para cierta constante $c_1$. Obviamente, si intercambiamos los parmetros $h_1,h_2,h_3$ apropiadamente, llegamos a la conclusin de que 
\[
\Delta_{h_1}^mf(z)=c_1,\  \Delta_{h_2}^mf(z)=c_2,\  \Delta_{h_3}^mf(z)=c_3
\]
para ciertas constantes $c_1,c_2,c_3$. Como $f$ es derivable, es fcil comprobar que
\[
\Delta_{h_1}^mf'(z)= \Delta_{h_2}^mf'(z)=\Delta_{h_3}^mf'(z)=0,
\]
y, ahora, aplicamos la hiptesis de induccin para concluir que $f'(z)$ es un polinomio algebraico de grado $\leq m-1$ o, lo que es lo mismo, que $f(z)$ es un polinomio de grado $\leq m$, que es lo que buscbamos.

{\hfill $\Box$}

Montel tambin estudi la ecuacin $\eqref{Pre_frepasofijo}$ para $X=\mathbb{R}^d$, con $d>1$,  y  $f:\mathbb{R}^d\to\mathbb{C}$ continua, y para $X=\mathbb{C}^d$ y $f:\mathbb{C}^d\to \mathbb{C}$ holomorfa.



\begin{teo}[Montel, varias variables] \label{Pre_montelvvcf} Supongamos que $h_1,\cdots,h_{\ell} \in \mathbb{R}^d$ son tales que $h_1\mathbb{Z}+\cdots+h_{\ell}\mathbb{Z}$ es un subconjunto denso de $\mathbb{R}^d$. 
Entonces:
\begin{itemize}
\item[$(i)$] Si $f\in C(\mathbb{R}^d,\mathbb{C})$ es tal que  $\Delta_{h_i}^m(f) =0$, $i=1,\cdots,\ell$. Entonces $f(x)=\sum_{|\alpha|<N}a_{\alpha}x^{\alpha}$ para cierto $N\in\mathbb{N}$, ciertos nmeros complejos $a_{\alpha}$, y todo $x\in\mathbb{R}^d$. Por tanto, $f$ es un polinomio algebraico con valores complejos en $d$ variables reales.

\item[$(ii)$] Supongamos que $d=2k$ e interpretamos los vectores $\{h_i\}_{i=1}^{\ell}$ como elementos de $\mathbb{C}^k=\mathbb{R}^{d}$. Si  $f:\mathbb{C}^k\to \mathbb{C}$ es holomorfa y  satisface las ecuaciones
$\Delta_{h_i}^m(f) =0$, $i=1,\cdots,\ell$, entonces  $f(z)=\sum_{|\alpha|<N}a_{\alpha}z^{\alpha}$ es un polinomio algebraico complejo, en $k$ variables complejas. 
\end{itemize}
\end{teo}

\begin{remark} Los subgrupos finitamente generados de  $(\mathbb{R}^d,+)$ que son densos en $\mathbb{R}^d$ se han estudiado en profundidad y, de hecho, existen varias caracterizaciones de los mismos. Por ejemplo, en \cite[Proposition 4.3]{W}, se puede consultar la demostracin del siguiente teorema:  
\begin{teo}\label{subgruposdensos} Sea $G=h_1\mathbb{Z}+h_2\mathbb{Z}+\cdots+h_{\ell}\mathbb{Z}$ el subgrupo aditivo de $\mathbb{R}^d$ generado por los vectores $\{h_1,\cdots,h_{\ell}\}$. Entonces las siguientes afirmaciones son equivalentes:
\begin{itemize}
\item[$(i)$] $G$ es un subconjunto denso de $\mathbb{R}^d$.
\item[$(ii)$] Si $\varphi:\mathbb{R}^d\to\mathbb{R}$ es una aplicacin lineal no nula, entonces $\varphi(G)\not\subseteq \mathbb{Z}$.
\item[$(iii)$] Si $\chi:\mathbb{R}^d\to\mathbb{R}/\mathbb{Z}$ es un homomorfismo continuo, entonces $\chi(G)\neq \{1\}$.
\item[$(iv)$] Si $h_k=(h_{1k},h_{2k},\cdots,h_{dk})$ son las coordenadas del vector $h_k$ respecto de la base cannica de $\mathbb{R}^d$ $(k=1,\cdots, \ell)$, entonces las matrices
\[
A(n_1,\cdots,n_\ell)=\left [
\begin{array}{cccccc}
h_{11} & h_{12} &  \cdots & h_{1\ell} \\
h_{21} & h_{22} & \cdots & h_{2\ell}\\
\vdots & \vdots & \ \ddots & \vdots \\
h_{d1} & h_{d2} &  \cdots &  h_{d\ell}\\
n_1 & n_2 & \cdots &  n_{\ell}
\end{array} \right].
\] 
son de rango $d+1$, para todo $(n_1,\cdots,n_\ell)\in\mathbb{Z}^\ell\setminus\{(0,\cdots,0)\}$.
\end{itemize}
\end{teo}

Un caso sencillo, que se ha considerado especialmente interesante, y que motiv histricamente el estudio de los subgrupos densos de $(\mathbb{R}^d,+)$, es el siguiente:
\end{remark}

\begin{cor}[Teorema de Kronecker] Dados $\theta_1,\theta_2,\cdots,\theta_d\in\mathbb{R}$, el grupo $\mathbb{Z}^d+(\theta_1,\theta_2,\cdots,\theta_d)\mathbb{Z}$ (que est generado por $d+1$ elementos) es un subconjunto denso de $\mathbb{R}^d$ si y slo si 
\[
n_1\theta_1+\cdots+n_d\theta_d\not\in \mathbb{Z}, \text{ para todo }(n_1,\cdots,n_d)\in\mathbb{Z}^{d}
\]  
Es decir, este grupo es denso en $\mathbb{R}^d$ si y slo si los nmeros $\{1,\theta_1,\cdots,\theta_d\}$ forman un sistema linealmente independiente sobre $\mathbb{Q}$.
\end{cor}

\noindent \textbf{Demostracin. } Los vectores $\{e_k\}_{k=1}^d\cup\{(\theta_1,\cdots,\theta_d)\}$, donde $$e_k=(0,0,\cdots,1^{(k\text{-sima posicin})},0,\cdots,0), \ \ k=1,\cdots, d,$$ generan el grupo $G=\mathbb{Z}^d+(\theta_1,\theta_2,\cdots,\theta_d)\mathbb{Z}$. Por tanto, el apartado $(iv)$ del Teorema \ref{subgruposdensos} nos dice que $G$ es denso en $\mathbb{R}^d$ si y slo si 
\begin{eqnarray*}
&\ & \det\left [
\begin{array}{cccccc}
\theta_{1} & 1 & 0 &  \cdots  & 0 \\
\theta_{2} & 0 & 1 &   \cdots  & 0 \\
\vdots & \vdots & \ \ddots & \vdots \\
\theta_{d} & 0 & 0&  \cdots  & 1 \\
n_0 & n_1 & n_2 &  \cdots &  n_{d}
\end{array} \right]  \\
&=&  (-1)^{d+2}n_0+\sum_{k=1}^d(-1)^{d+2+k}n_k(-1)^{k+1}\theta_k \det(I_{d-1}) \\
  &=&  (-1)^{d}\left(n_0-\sum_{k=1}^dn_k\theta_k\right) \neq 0
\end{eqnarray*}
para todo $(n_0,n_1,\cdots,n_d)\in\mathbb{Z}^{d+1}$, lo cual es obviamente equivalente a que  
\[
n_1\theta_1+\cdots+n_d\theta_d\not\in \mathbb{Z}, \text{ para todo }(n_1,\cdots,n_d)\in\mathbb{Z}^{d},
\]  
que es lo que buscbamos. Para demostrar la ltima afirmacin del teorema, basta observar que, si $\{1,\theta_1,\cdots,\theta_d\}$ forma un sistema linealmente dependiente sobre $\mathbb{Q}$, entonces existen nmeros racionales $r_i=n_i/m_i$, $i=0,1,\cdots, d$, tales que 
\[
r_0+r_1\theta_1+\cdots+r_d\theta_d= 0, 
\]  
y, por tanto, si multiplicamos por $m=\prod_{k=0}^dm_k$ a ambos lados de la ecuacin, obtendremos que
\[
n_0^*+n_1^*\theta_1+\cdots+n_d^*\theta_d= 0
\]
para ciertos nmeros enteros $n_0^*,n_1^*,\cdots,n_d^*$. En particular, 
\[
n_1^*\theta_1+\cdots+n_d^*\theta_d\in\mathbb{Z}.
\]  
{\hfill $\Box$}

En general, se sabe (ver \cite[Theorem 3.1]{W})  que si $G$ es un subgrupo finitamente generado de $(\mathbb{R}^d,+)$, entonces la clausura topolgica de $G$ admite una descomposicin del tipo 
$$\overline{G}^{\mathbb{R}^d}=V\oplus \Lambda,$$
donde $V$ es un subespacio vectorial de $\mathbb{R}^d$ y $\Lambda$ es un subgrupo discreto de $\mathbb{R}^d$.

\begin{remark} Los Teoremas de Montel en una variable (Teorema $\ref{Pre_montel_1v}$ ) y en varias variables (Teorema $\ref{Pre_montelvvcf}$)  han sido recientemente generalizados por Almira \cite{A_NFAO} y Almira y Abu-Helaiel \cite{AK_CJM} al caso de distribuciones y tambin al caso discreto (cuando las funciones estn definidas en $\mathbb{Z}^{n}$) utilizando tcnicas de anlisis funcional. Finalmente, en \cite{AS_montel}, Almira y Szkelyhidi han demostrado, con tcnicas distintas, un Teorema tipo Montel para funciones definidas en grupos Abelianos finitamente generados y, adems, han mejorado el Teorema de Montel en varias variables al demostrar que se verifica el siguiente resultado:

\begin{teo} 
Sea $t$ un entero positivo, sean $n_1,n_2,\dots,n_t$ nmeros naturales, y, adems, sea $f:\mathbb{R}^d\to\mathbb{R}$ una funcin continua que satisface  
\[
\Delta_{h_k}^{n_k+1}f(x)=0
\]
para todo $x$ en $\mathbb{R}^d$ y para $k=1,\cdots,t$. Si el subgrupo aditivo $G$  de $\mathbb{R}^d$ generado por $\{h_1,h_2,\dots,h_t\}$ satisface $\overline{G}=V\oplus \Lambda$, donde $V$ es un subespacio vectorial de $\mathbb{R}^d$, y  $\Lambda$ es un subgrupo aditivo discreto de  $\mathbb{R}^d$, entonces existen polinomios algebraicos ordinarios  $p_{\lambda}:\mathbb{R}^d\to\mathbb{R}$, con $\lambda$ en $\Lambda$, tales que
\[
f(x+\lambda)=p_{\lambda}(x) \text{ para todo } x\in V \text{ y todo } \lambda \in\Lambda. 
\]
 Adems, se tiene que el grado total de $p_{\lambda}$ satisface la desigualdad
\begin{equation*}
\deg p_{\lambda}\leq  n_1+n_2+\dots+n_t+t-1
\end{equation*}
para todo $\lambda$ en $\Lambda$. En particular, si $V=\mathbb{R}^d$, entonces $f$ es un polinomio algebraico ordinario.  Adems, si $d=1$ y $V=\mathbb{R}$, el grado de $f$ es menor o igual que  $\min\{n_k:k=1,\cdots,t\}$. 
\end{teo}

\end{remark}

\section{El teorema de Montel-Popoviciu}

Aunque el artculo de Montel \cite{montel} fue publicado en 1937, l ya haba obtenido sus resultados en 1935 y, de hecho, ese ao imparti un seminario en el que se explicaba su teorema, en el departamento de matemticas de la Universidad Politcnica de Cluj Napoca, en Rumana. En dicho seminario estaba presente su alumno de doctorado, el matemtico rumano Tiberiu Popoviciu, quien haba defendido su tesis doctoral en la Escuela Normal Superior de Pars en 1933  \cite{popoviciu_tesis}.  Popoviciu capt inmediatamente las ideas de Montel y, de hecho, en 1936 public un artculo \cite{popoviciu1} en el que se mejoraba sensiblemente el resultado original. Concretamente, demostr el siguiente teorema:  

\begin{teo}[Montel-Popoviciu, 1935] \label{MP_MP} Sea $f:\mathbb{R}\to\mathbb{R}$ una solucin del sistema de ecuaciones funcionales
\[
\Delta^{m+1}_{h_1}f(x)=\Delta^{m+1}_{h_2}f(x)=0,  \text{  } (x\in\mathbb{R}),
\]
donde $h_1,h_2\in\mathbb{R}\setminus \{0\}$. Entonces:
\begin{itemize}
\item[(i)] Para cada $x_0\in\mathbb{R}$ existe un nico polinomio $P(x,y)\in \Pi_m^2$ tal que $f(x_0+ih_1+jh_2)=P(ih_1,jh_2)$ para todo $(i,j)\in\mathbb{Z}^2$.
\item[(ii)] Si $h_1/h_2\not\in\mathbb{Q}$ y existe un intervalo abierto no vaco $I=(a,b)$ tal que $f_{|I}$ es una funcin acotada, entonces el polinomio $P(x,y)$ que aparece en $(i)$ es de la forma $P(x,y)=A_0(x+y)$ para cierto polinomio $A_0(t)$ de grado $\leq m$ y, consecuentemente, 
\[
\Delta^{m+1}_{ih_1+jh_2}f(x)=0,  \text{ para todo } x\in\mathbb{R} \text{, y todo } i,j\in\mathbb{Z}.  
\]
\item[(iii)]  Si $h_1/h_2\not\in\mathbb{Q}$ y $f$ es continua en al menos $m+1$ puntos, entonces $f(x)=A_{0}(x)$ $(x\in\mathbb{R})$, es un polinomio algebraico de grado $\leq m$.
\end{itemize}
\end{teo}

Vamos a seguir el espritu de la prueba original de Popoviciu, aunque nosotros explicaremos ms detalladamente todas nuestras construcciones, algunas de las cuales difieren de forma significativa de lo que se afirma -a veces, sin demostracin- en el artculo original.  Esta prueba servir, posteriormente, para obtener una interesante descripcin cualitativa del grafo de los polinomios discontinuos. 

Comenzamos, pues, con algunos resultados tcnicos sobre polinomios de una y dos variables reales. El primero es un resultado estndar sobre localizacin de ceros.

\begin{lem} \label{MP_ceros_pol}
Sea $p(z)=a_0+a_1z+\cdots+a_Nz^N\in \mathbb{C}[z]$ un polinomio algebraico de grado exactamente $N$ (i.e., $a_N\neq 0$)  con coeficientes complejos, y sea $\xi\in\mathbb{C}$ cualquiera de sus ceros. Entonces
\[
|\xi|\leq \max\{1,\sum_{k=0}^{N-1}\frac{|a_k|}{|a_N|}\}.
\] 
\end{lem}

\noindent \textbf{Demostracin. } Si $|\xi|\leq 1$, ya hemos acabado. Supongamos, pues, que $|\xi|>1$. Como $p(\xi)=0$, tenemos que tambin $\frac{1}{|a_N|}p(\xi)=0$, de modo que, si despejamos el trmino lder de esta expresin, y tomamos valor absoluto a ambos lados de la igualdad, obtenemos que
\begin{eqnarray*}
|\xi|^N &=& | \frac{a_0}{a_N}+\frac{a_1}{a_N}\xi+\cdots+\frac{a_{N-1}}{a_N}\xi^{N-1}| \leq \sum_{k=0}^{N-1}\frac{|a_k|}{|a_N|}\max\{1,|\xi|,\cdots,|\xi|^{N-1}\}\\
&=& \sum_{k=0}^{N-1}\frac{|a_k|}{|a_N|}|\xi|^{N-1},
\end{eqnarray*}
por tanto, $|\xi|\leq  \sum_{k=0}^{N-1}\frac{|a_k|}{|a_N|}$, que es lo que buscbamos. {\hfill $\Box$}

A continuacin demostramos un resultado de naturaleza ms tcnica.

\begin{lem}\label{MP_poli_banda}
Sea $p(z)=a_0+a_1z+\cdots+a_Nz^N\in \mathbb{C}[z]$ un polinomio algebraico de grado exactamente $N$ (i.e., $a_N\neq 0$)  con coeficientes complejos, y supongamos que $N\geq 1$.  Sea adems $\{q_n(z)\}_{n=1}^\infty$ es una sucesin de polinomios de grado $\leq N$ tales que $q_n(z)=a_{0n}+a_{1n}z+\cdots+a_{Nn}z^N$, y $\max\{|a_k-a_{kn}|:k=0,1,\cdots,N\}<|a_N|/2$, $n=1,2,\cdots,\infty$. Entonces, si $\{x_n\}$ es una sucesin de nmeros tal que $|x_n|\to+\infty$, entonces $|q_n(x_n)|\to\infty$. 
\end{lem}
\noindent \textbf{Demostracin. } Sea $n\in\mathbb{N}$ y tomemos $\xi$ un cero de $q_n(z)$. Gracias al Lema \ref{MP_ceros_pol}, sabemos que
\[
|\xi|\leq \max\{1,\sum_{k=0}^{N-1}\frac{|a_{kn}|}{|a_{Nn}|}\},
\]
y, como $|a_{kn}|\leq |a_{kn}-a_k|+|a_k|\leq \frac{|a_N|}{2}+|a_k|$, $|a_{Nn}|\geq \frac{a_N}{2}$, tenemos que 
\[
|\xi|\leq \max\{1,\sum_{k=0}^{N-1}\frac{2(\frac{|a_N|}{2}+|a_k|)}{|a_{N}|}\}=:M.
\]
Es decir, todos los ceros de $q_n(z)$ (para todo $n$) estn en $B_M=\{z\in\mathbb{C}:|z|\leq M\}$. Si $|x_n|\to\infty$, entonces obviamente $\mathbf{\text{dist}}(x_n,B_M)\to\infty$. Ahora, si $\{\alpha_{kn}\}_{k=1}^N$ denota el conjunto de los ceros de $q_n(z)$, entonces obviamente $q_n(z)=a_{Nn}\prod_{k=1}^N(z-\alpha_{kn})$ y, por tanto,
\[
|q_n(x_n)|=|a_{Nn}|\prod_{k=1}^N|x_n-\alpha_{kn}|\geq \frac{|a_N|}{2}(\text{dist}(x_n,B_M))^N\to\infty. \ \ (n\to\infty).
\]
{\hfill $\Box$}

\begin{cor}\label{MP_pol_2v}
Sea $P(x,y)\in\Pi_{m}^2$ un polinomio algebraico de grado menor o igual que $m$ en ambas variables. Supongamos que existen $m+1$ nmeros distintos $\{\alpha_k\}_{k=1}^{m+1}$ y sucesiones  de puntos del plano $\{(u_{k,n},v_{k,n})\}$ tales que $u_{k,n}+v_{k,n}\to \alpha_k$ para $n\to\infty$ para $k=1,2,\cdots,m+1$. Supongamos, adems, que $|u_{k,n}|\to\infty$ para $n\to\infty$ y que $\{P(u_{k,n},v_{k,n})\}_{n=1}^\infty$ es una sucesin acotada para cada $k$. Entonces 
\[
P(x,y)=A_0(x+y),
\]
donde $A_0$ es un polinomio en una variable, de grado $\leq m$.
\end{cor}

\noindent \textbf{Demostracin. } Consideremos el cambio de variables $\varphi(x,y)=(x,x+y)$. Si denotamos $f_1=x$, $f_2=x+y$, entonces $y=f_2-f_1$ y, consecuentemente, todo polinomio $P(x,y)=\sum_{i=0}^m\sum_{j=0}^ma_{i,j}x^iy^j\in \Pi_m^2$ admite una representacin del tipo
\begin{eqnarray*}
P(x,y) &= & \sum_{i=0}^m\sum_{j=0}^ma_{i,j}x^iy^j = \sum_{i=0}^m\sum_{j=0}^ma_{i,j}f_1^i(f_2-f_1)^j \\
&=& \sum_{i=0}^m\sum_{j=0}^ma_{i,j}f_1^i (\sum_{s=0}^j\binom{j}{s}(-1)^{j-s}f_2^sf_1^{j-s}) \\
&=& \sum_{i=0}^m\sum_{j=0}^m(\sum_{s=0}^ja_{i,j}\binom{j}{s}(-1)^{j-s}f_2^sf_1^{i+j-s}) \\
&=& \sum_{i=0}^{2m}A_i(f_2)f_1^i \\
&=& \sum_{i=0}^{2m}A_i(x+y)x^i,
\end{eqnarray*}
donde $A_i(t)$ es un polinomio en una nica variable, de grado $\leq m$, para $i=0,1,\cdots, 2m$. Una vez hemos representado el polinomio $P(x,y)$ de la forma anterior, tomamos $\{(u_{k,n},v_{k,n})\}$ y $\{\alpha_k\}$ verificando las hiptesis del corolario. Sea $N=\max\{i\in \{0,\cdots, 2m\}:A_i\neq 0\}$. Entonces  
\[
P(x,y)= \sum_{i=0}^{N}A_i(x+y)x^i \text{ y } A_N\neq 0.
\]
Queremos demostrar que $N=0$. Supongamos, por el contrario, que $N>1$. Vamos a demostrar que $A_N(\alpha_k)=0$ para $k=1,2,\cdots, m+1$, lo cual implicara que $A_N=0$ pues $A_N(t)$ es un polinomio de grado $\leq m$. Supongamos, por el contrario, que $A_N(\alpha_k)\neq 0$ para algn $k$. Consideremos el polinomio 
\[
p(z)=\sum_{i=0}^{N}A_i(\alpha_k)z^i 
\]
Como $\{u_{k,n}+v_{k,n}\}\to \alpha_k$ para 
$n\to \infty$, y las funciones $A_i(t)$ son continuas, sabemos que $\{A_i(u_{k,n}+v_{k,n})\}\to A_i(\alpha_k)$ para todo $i\in \{0,1,\cdots, N\}$. Es ms, podemos asumir sin prdida de generalidad que $|A_i(u_{k,n}+v_{k,n})-A_i(\alpha_k)|<|A_N(\alpha_k)|/2$, $i=0,1,\cdots,N$, $n\in\mathbb{N}$. Si consideramos los polinomios
\[
q_n(z)=\sum_{i=0}^{N}A_i(u_{k,n}+v_{k,n})z^i, \ n=1,2,\cdots,
\]
y usamos el Lema \ref{MP_poli_banda}, obtenemos que $|q_n(u_{k,n})|\to\infty$ para $n\to\infty$. Sin embargo, $|q_n(u_{k,n})|=|P(u_{k,n},v_{k,n})|$, que es una sucesin acotada. Esto nos lleva a la conclusin de que $A_N(\alpha_k)=0$ forzosamente, que es lo que buscbamos. {\hfill $\Box$}

\noindent \textbf{Demostracin del Teorema de Montel-Popoviciu (Teorema \ref{MP_MP}). }

\noindent $(i)$. Sea $f:\mathbb{R}\to\mathbb{R}$ tal que 
\[
\Delta_{h_1}^{m+1}f(x)=\Delta_{h_2}^{m+1}f(x)=0 \ \ \text{ para todo } x\in\mathbb{R},
\]
Sea $x_0\in\mathbb{R}$ y consideremos el polinomio de interpolacin de Lagrange (en dos variables $x,y$) que interpola la tabla de valores $f_{i,j}=f(x_0+ih_1+jh_2)$, $i,j=0,1,\cdots, m$, en los nodos $(ih_1,jh_2)$, $i,j=0,1,\cdots, m$. Es decir, 
\[
P(ih_1,jh_2)=f(x_0+ih_1+jh_2), \ i,j=0,1,\cdots, m.
\]
Vamos a demostrar que $f(x_0+ih_1+jh_2)=P(ih_1,jh_2)$ para todo $i,j\in\mathbb{Z}$.  Para verlo, observemos que
\begin{eqnarray*}
0 &=& \Delta_{h_1}^{m+1}f(x_0+jh_2)=\sum_{k=0}^{m+1}\binom{m+1}{k}(-1)^{m+1-k}f(x_0+kh_1+jh_2) \\
&=& \sum_{k=0}^{m}\binom{m+1}{k}(-1)^{m+1-k}f(x_0+kh_1+jh_2) +f(x_0+(m+1)h_1+jh_2)\\
\end{eqnarray*}
Por otra parte, como $P(x,y)=\sum_{t=0}^m\sum_{s=0}^m\alpha_{t,s}x^ty^s$, si fijamos $y=jh_2$ y consideramos la funcin $g(x)=P(x,jh_2)$, resulta que $g$ es un polinomio de grado $\leq m$, por lo que $\Delta_h^{m+1}g(x)=0$ para todos los valores $x,h\in\mathbb{R}$, de modo que  
\begin{eqnarray*}
0 &=& \Delta_{h_1}^{m+1}g(0)=\sum_{k=0}^{m+1}\binom{m+1}{k}(-1)^{m+1-k}g(kh_1) \\
&=& \sum_{k=0}^{m}\binom{m+1}{k}(-1)^{m+1-k}P(kh_1,jh_2) +P((m+1)h_1,jh_2)\\
&=& \sum_{k=0}^{m}\binom{m+1}{k}(-1)^{m+1-k}f(x_0+kh_1+jh_2) +P((m+1)h_1,jh_2).
\end{eqnarray*}
Se sigue que $P((m+1)h_1,jh_2)=f(x_0+(m+1)h_1+jh_2)$. Si hubiramos partido del polinomio de Lagrange que interpola los valores $f_{i,j}=f(x_0+ih_1+jh_2)$ para $i=1,2,\cdots,m+1, j=0,1,\cdots,m$ en los nodos $(ih_1,jh_2)$ (que, como acabamos de demostrar, coincide con nuestro polinomio $P(x,y)$), el mismo argumento nos llevara a la conclusin de que tambin $P((m+2)h_1,jh_2)=f(x_0+(m+2)h_1+jh_2)$. De forma similar, despejando esta vez el primer trmino de la sumas en vez del ltimo, y tomando como punto de partida la igualdad $\Delta_{h_1}^{m+1}f(x_0 
-h_1+jh_2)=0$, podramos concluir que $P(-h_1,jh_2)=f(x_0-h_1+jh_2)$. Repitiendo estos argumentos infinitas veces (tanto hacia delante como hacia atrs), se obtiene que $P(kh_1,jh_2)=f(x_0+kh_1+jh_2)$ para todo $k\in\mathbb{Z}$ y $j=0,1,\cdots,m$. Ahora, tomando la funcin $h(y)=P(ih_1,y)$, podemos argumentar de modo similar que, en realidad, 
\[
P(ih_1,jh_2)=f(x_0+ih_1+jh_2) \text{ para todo } (i,j)\in\mathbb{Z}^2.
\] 
Por otra parte, el polinomio $P(x,y)=\sum_{t=0}^m\sum_{s=0}^m\alpha_{t,s}x^ty^s$ para ciertas constantes $\alpha_{t,s}$, por lo que, gracias a lo ya expuesto en la demostracin del Corolario \ref{MP_pol_2v}, tambin admite una expresin del tipo 
\[
P(x,y)=\sum_{i=0}^{2m}A_i(x+y)x^i,
\]
con $A_i(t)$ un polinomio de grado $\leq m$ para cada $i$. 



\noindent $(ii)$. Supongamos ahora que $h_1/h_2\not\in\mathbb{Q}$ y $f$ est acotada en un cierto intervalo $(a,b)$ con $a<b$. Como $h_1/h_2$ es irracional, sabemos que los puntos $\{x_0+ih_1+jh_2\}_{i,j\in\mathbb{Z}}$ forman un subconjunto denso de la recta y, por tanto, fijados $m+1$ elementos distintos $\{\alpha_k\}_{k=1}^{m+1}\subset (a,b)$,  existen infinitos valores enteros $i_{k,n},j_{k,n}$ tales que $\bigcup_{k\leq m+1}\bigcup_{n\in\mathbb{N}}\{x_0+i_{k,n}h_1+j_{k,n}h_2\} \subseteq (a,b)$, $\lim_{n\to\infty}(i_{k,n}h_1+j_{k,n}h_2)=\alpha_k-x_0$, $\lim_{n\to\infty}|i_{k,n}h_1|=\infty$, $k=1,2,\cdots, m+1$. Como $P(ih_1,jh_2)=f(x_0+ih_1+jh_2)$  para todo $(i,j)\in\mathbb{Z}^2$, podemos utilizar el Corolario \ref{MP_pol_2v}, para concluir que $P(x,y)=A_0(x+y)$. En particular, 
\[
f(x_0+ih_1+jh_2)=A_0(ih_1+jh_2) \text{ para todo } (i,j)\in\mathbb{Z}^2.
\]
Esto, unido a que $A_0(t)$ es un polinomio de grado $\leq m$, y a que todos los clculos anteriores son vlidos independientemente del punto $x_0$, implica que 
\[
\Delta_{ih_1+jh_2}^{m+1}f(x)=0 \text{ para todo } x\in\mathbb{R} \text{ y todo } (i,j)\in\mathbb{Z}^2.
\]

\noindent $(iii)$. Supongamos ahora que  $h_1/h_2\not\in\mathbb{Q}$ y  $f$ es continua en $m+1$ puntos distintos $\{s_k\}_{k=1}^{m+1}$. Esto obviamente implica que $f$ est acotada en un intervalo abierto no vaco y, por tanto, para cada $x_0\in\mathbb{R}$, existe un polinomio $A_{x_0}(t)$ de grado $\leq m$ tal que $f(x_0+ih_1+jh_2)=A_{x_0}(ih_1+jh_2) $ para todo $(i,j)\in\mathbb{Z}^2$. Decir que $f$ es un polinomio algebraico de grado $\leq m$ equivale, por tanto, a decir que $A_{x_0}(x+x_1-x_0)=A_{x_1}(x)$, independientemente de los valores de $x_0$ y $x_1$. Para ver esto, basta observar que, al ser $x_1$ arbitrario, y teniendo en cuenta la definicin de $A_{x_1}$, sabemos que
\[
f(x_1)=A_{x_1}(0)=A_{x_0}(x_1-x_0),
\]
es decir, $f(t)=A_{x_0}(t-x_0)$ para todo valor $t\in\mathbb{R}$, que es lo que buscamos.

Ahora bien, hay que probar la identidad $A_{x_0}(x+x_1-x_0)=A_{x_1}(x)$ o, lo que es equivalente, $A_{x_0}(x)=A_{x_1}(x+x_0-x_1)$.  Tomemos $i_{k,n},j_{k,n},i_{k,n}^*,j_{k,n}^* \in\mathbb{Z}$ tales que $x_0+i_{k,n}h_1+j_{k,n}h_2\to s_k$, y  $x_1+i_{k,n}^*h_1+j_{k,n}^*h_2\to s_k$, para $k=1,2,\cdots, m+1$, y consideremos el polinomio $C(x)=A_{x_1}(x+x_0-x_1)$. 
Entonces 
\begin{eqnarray*}
f(s_k) &=& \lim_{n\to\infty}f(x_0+i_{k,n}h_1+j_{k,n}h_2)\\
&=& \lim_{n\to\infty}A_{x_0}(i_{k,n}h_1+j_{k,n}h_2)=A_{x_0}(s_k-x_0),\ 1\leq k\leq m+1.
\end{eqnarray*}
Por otra parte,
\begin{eqnarray*}
f(s_k)&=& \lim_{n\to\infty}f(x_1+i_{k,n}^*h_1+j_{k,n}^*h_2)\\
&=& \lim_{n\to\infty}A_{x_1}(i_{k,n}^*h_1+j_{k,n}^*h_2)=A_{x_1}(s_k-x_1),\ 1\leq k\leq m+1.
\end{eqnarray*}
Por tanto, 
\[
C(s_k-x_0)=A_{x_1}(s_k-x_0+x_0-x_1)=A_{x_1}(s_k-x_1)=A_{x_0}(s_k-x_0),\ 1\leq k\leq m+1.
\]
Es decir, $C(x)$ y $A_{x_0}$ coinciden en al menos $m+1$ puntos, por lo que son en realidad el mismo polinomio. Por tanto, $A_{x_0}(x)=A_{x_1}(x+x_0-x_1)$, que es lo que queramos demostrar. 

{\hfill $\Box$}

\section{Descripcin cualitativa de los grafos de los polinomios discontinuos}

La tcnica de interpolacin utilizada por Popoviciu para su prueba del Teorema de Montel (mejorado), se puede aprovechar para el estudio cualitativo del grafo de un polinomio discontinuo.  Concretamente, se puede demostrar el siguiente resultado, que generaliza el Teorema de Darboux para la ecuacin de Cauchy al contexto de los polinomios generalizados (el contenido de esta seccin se basa muy fuertemente en nuestro artculo \cite{AK_D}, pero en el caso de funciones reales de una variable real, se pueden consultar tambin los trabajos anteriores \cite{almira_antonio} y \cite{AK_MJM}):

\begin{teo}[Descripcin de $G(f)$ para funciones de una variable]  \label{MP_grafo_1v} Si $f:\mathbb{R}\to\mathbb{R}$ satisface la ecuacin funcional de Frchet
\[
\Delta_h^{m+1}f(x)=0 \ \ \text{ para todo }(x,h) \in\mathbb{R}^2,
\]
y no es un polinomio ordinario, entonces   $f$ no est acotada en ningn intervalo abierto no vaco. Adems, para todo  $x\in\mathbb{R}$ existe un intervalo no acotado $I_x\subseteq \mathbb{R}$ tal que $\{x\}\times I_x\subseteq  \overline{G(f)}^{\mathbb{R}^2}$. Finalmente, $ \overline{G(f)}^{\mathbb{R}^2}$ contiene un abierto no acotado. 
\end{teo}

\noindent \textbf{Demostracin. }  Si $f:\mathbb{R}\to\mathbb{R}$ es una solucin de la ecuacin de Frchet $\Delta_h^{m+1}f=0$, entonces el apartado $(i)$ del Teorema \ref{MP_MP} se puede mejorar, al demostrar que, en realidad,   
\begin{itemize}
\item[$(i)^*$] Para cada $x_0\in\mathbb{R}$, y para cada $h_1,h_2\in\mathbb{R}$, existe un nico polinomio $P(x,y)\in \Pi_m^2$ tal que $f(x_0+rh_1+sh_2)=P(rh_1,sh_2)$ para todo $(r,s)\in\mathbb{Q}^2$.
\end{itemize}
Para verlo, basta comprobar que, si aplicamos $(i)$ del Teorema \ref{MP_MP} a $f$ cambiando los pasos $h_1,h_2$ por $h_1^*=h_1/p, h_2^*=h_2/q$, respectivamente, donde $p,q$ son dos nmeros enteros (no nulos) arbitrarios, el nuevo polinomio, que denotamos por $P^*$, coincide con el polinomio original, $P(x,y)$, en una malla infinita de puntos, pues 
\[
P(ih_1,jh_2)=f(x_0+ih_1+jh_2)=f(x_0+iph_1^*+jqh_2^*)=P^*(iph_1^*,jqh_2^*)=P^*(ih_1,jh_2) 
\]
para todo $i,j\in\mathbb{Z}$, por lo que en realidad son el mismo polinomio, lo cual implica obviamente $(i)^*$, pues $p,q$ eran enteros (no nulos) arbitrarios.

Se sigue que, si $\Delta^{m+1}_hf=0$ para todo $h$ y $P_{x_0,h_1,h_2}(x,y)$ denota al polinomio de $\Pi_{m}^2 $ que satisface $(i)^*$, entonces
 \begin{equation}\label{MP_grafo_1}
 \Gamma_{x_0,h_1,h_2}:=\{(x_0+u+v,P_{x_0,h_1,h_2}(u,v)):u,v\in\mathbb{R}\}\subseteq \overline{G(f)}^{\mathbb{R}^2}, 
 \end{equation}
pues $\mathbb{Q}$ es un subconjunto denso de $\mathbb{R}$. Lo interesante ahora es, pues, estudiar los conjuntos $\Gamma_{x_0,h_1,h_2}$. Si $P_{x_0,h_1,h_2}(x,y)=A(x+y)$ para cierto polinomio $A$, entonces el conjunto $\Gamma_{x_0,h_1,h_2}$ tiene interior vaco y, de hecho, no es ms que el grafo de un polinomio algebraico ordinario, por lo que, en este caso, la propiedad  (\ref{MP_grafo_1}) no dice gran cosa. Sin embargo, podemos demostrar que, si $f$ es solucin de la ecuacin de Frchet $\Delta^{m+1}_hf=0$ y no es un polinomio algebraico, entonces existen $x_0,h_1,h_2\in\mathbb{R}$ tales que $P_{x_0,h_1,h_2}(x,y)$ no es un polinomio en la variable $x+y$. Concretamente, demostramos que, para valores apropiados de  $x_0,h_1$ y $h_2$, 
\begin{equation}\label{MP_poli_estricto}
P_{x_0,h_1,h_2}(x,y)=\sum_{i=0}^NA_i(x+y)x^i,\ \text{ con } A_N(t)\neq 0 \text{ y } N\geq 1,
\end{equation}
donde $A_i(t)$ es un polinomio de grado $\leq m$ para $i=0,1,\cdots,N$.  Evidentemente, si se satisface $\eqref{MP_poli_estricto}$, entonces, para cada $\alpha\in\mathbb{R}\setminus Z(A_N)$ (donde $Z(A_N)=\{s\in\mathbb{R}:A_N(s)=0\}$ es un conjunto de a lo sumo $m$ puntos), tendremos que $p_{\alpha}(x)=P_{x_0,h_1,h_2}(x,\alpha-x)=\sum_{i=0}^NA_i(\alpha)x^i$ es un polinomio no constante y, por tanto, $p_{\alpha}(\mathbb{R})$ es un intervalo no acotado. Adems, 
\[
\{x_0+\alpha\}\times p_{\alpha}(\mathbb{R})\subseteq \Gamma_{x_0,h_1,h_2}\subseteq \overline{G(f)}^{\mathbb{R}^2}.
\]
Es decir, si demostramos $\eqref{MP_poli_estricto}$, entonces $f$ no podr ser localmente acotada y, para todo  $x\in\mathbb{R}$ existir un intervalo no acotado $I_x\subseteq \mathbb{R}$ tal que $\{x\}\times I_x\subseteq  \overline{G(f)}^{\mathbb{R}^2}$. 

Veamos, por ltimo, que si $P=P_{x_0,h_1,h_2}$ satisface  $\eqref{MP_poli_estricto}$, entonces $\Gamma_{x_0,h_1,h_2}$ contiene un abierto no acotado. Para comprobarlo, definimos la funcin $\varphi:\mathbb{R}^2\to\mathbb{R}^2$, $$\varphi(x,y)=(x+y+x_0,P(x,y)).$$ Un sencillo clculo nos revela que $$\det \varphi'(x,y)=P_y-P_x=-\sum_{k=1}^NkA_k(x+y)x^{k-1}$$ es un polinomio no nulo y, por tanto, $\Omega=\mathbb{R}^2\setminus \{(x,y): \det \varphi'(x,y)=0\}$ es un abierto no vaco del plano (de hecho, es un abierto denso en el plano). Podemos, pues, utilizar el teorema de la aplicacin abierta para funciones diferenciales con la funcin $\varphi$, concluyendo que $W=\varphi(\Omega)$ es un abierto de $\mathbb{R}^2$ contenido en $\Gamma_{x_0,h_1,h_2}$. Las inclusiones $\{x_0+\alpha\}\times p_{\alpha}(\mathbb{R})\subseteq \Gamma_{x_0,h_1,h_2}$ demuestran, adems, que $W$ no es un conjunto acotado. 

Veamos ahora que, en efecto, se satisface la relacin $\eqref{MP_poli_estricto}$ para cierta eleccin de $x_0,h_1,h_2$. Si esto no fuera as, entonces, fijados $h_1,h_2$, todos los polinomios  $P_{x_0,h_1,h_2}(x,y)$ seran de la forma $P_{x_0,h_1,h_2}(x,y)=A_{x_0}(x+y)$ para cierto polinomio $A_{x_0}$ de grado $\leq m$. Como hemos supuesto que $f$ no es un polinomio algebraico, esto implicar la existencia de dos puntos distintos $x_0,x_1$ tales que $A_{x_1}(0)\neq A_{x_0}(x_1-x_0)$.  Consideremos ahora el polinomio $P_{x_0,x_1-x_0,h_2}(x,y)$, el cual, por hiptesis, debe ser de la forma $P_{x_0,x_1-x_0,h_2}(x,y)=A(x+y)$ para cierto polinomio $A(t)$ de grado $\leq m$. Un simple clculo nos conduce a que 
\[
A(x_1-x_0)=P_{x_0,x_1-x_0,h_2}(x_1-x_0,0)=f(x_0+(x_1-x_0))=f(x_1)=A_{x_1}(0).
\]
Por otra parte, para cada $j\in\mathbb{Z}$, tenemos que 
\[
A(jh_2)=P_{x_0,x_1-x_0,h_2}(0,jh_2)=f(x_0+jh_2)=A_{x_0}(jh_2),
\]
por lo que $A$ y $A_{x_0}$ coinciden en infinitos puntos y son, en consecuencia, el mismo polinomio. Se sigue que $A_{x_1}(0)= A_{x_0}(x_1-x_0)$, lo cual contradice la hiptesis de que $f$ no es un polinomio.

{\hfill $\Box$}

El resultado anterior puede generalizarse al caso de varias variables:

\begin{teo}[Descripcin de $G(f)$ para funciones de varias variables]  \label{MP_grafo_vv} Si $f:\mathbb{R}^n\to\mathbb{R}$ satisface la ecuacin funcional de Frchet
\[
\Delta_h^{m+1}f(x)=0 \ \ \text{ para todo }x,h \in\mathbb{R}^n,
\]
y $f(x_1,\cdots,x_n)$ no es un polinomio ordinario, entonces   $f$ no est acotada en ningn abierto no vaco. Adems,  $\overline{G(f)}^{\mathbb{R}^{n+1}}$ contiene un abierto no acotado. 
\end{teo}

Antes de abordar la prueba del Teorema \ref{MP_grafo_vv}, es necesario introducir algunas construcciones, as como varios resultados tcnicos relacionadas con ellas. Para empezar, observemos que si $f$ no es un polinomio algebraico en sentido ordinario, entonces existen $s\in \{1,\cdots,n\}$ y  $(a_1,a_2,\cdots,a_{s-1},a_{s+1},\cdots,a_n)\in\mathbb{R}^{n-1}$  (valores que dejamos fijos a partir de ahora) tales que $$g(x)=f(a_1,a_2,\cdots,a_{s-1},x,a_{s+1},\cdots,a_n)$$ no es un polinomio algebraico. Esto ha sido demostrado de varias formas y puede consultarse, por ejemplo, en   \cite{AK_CJM},  \cite{kuczma}, \cite{prager}. Adems, si tenemos en cuenta la demostracin del Teorema \ref{MP_grafo_1v}, sabemos que, si denotamos por $p_{x_0,\alpha,\beta}(x,y)$ al nico polinomio de 
$\Pi_{m,\max}^2$ tal que 
\[
p_{x_0,\alpha,\beta}(i\alpha,j\beta)=g(x_0+i\alpha+j\beta),\ \ \text{ para } i,j=0,1,\cdots,m,
\]
entonces existen $a_s,h_s,h_{n+1}\in\mathbb{R}$, $1\leq N\leq 2m$, y polinomios $A_k\in \Pi_m$, $k=0,1,\cdots, N$ tales que 
\[
p_{a_s,h_s,h_{n+1}}(x,y)=\sum_{k=0}^NA_k(x+y)x^k,\ \ \text{ donde } A_N\neq 0.
\]
Dejamos tambin fijos en todo lo que sigue los valores  $a_s,h_s,h_{n+1}$, y tomamos $a=(a_1,\cdots,a_n)\in\mathbb{R}^n$, $h_1,\cdots,h_{s-1},h_{s+1}, h_{n+1}\in\mathbb{R}\setminus\{0\}$ (tambien fijos), y $\gamma=\{v_k\}_{k=1}^{n+1}\subset \mathbb{R}^n\setminus \{(0,0,\cdots,0)\}$. Entonces, gracias a la tcnica del producto tensorial, es fcil demostrar que existe un nico polinomio algebraico $P(t_1,\cdots,t_{n+1})\in \Pi_{m,\max}^{n+1}$ tal que
\[
P(i_1h_1,i_2h_2,\cdots,i_{n+1}h_{n+1})=f(a+\sum_{k=1}^{n+1}i_kh_kv_k), \text{ para } 0\leq i_k\leq m,\ 1\leq k\leq n+1.
\] 
Dicho polinomio lo denotamos, a partir de ahora, por $P_{a,h,\gamma}$, donde $h:=(h_1,\cdots,h_{n+1})$.  La demostracin del Teorema \ref{MP_grafo_vv} se apoya en varios lemas que enunciamos y demostramos a continuacin. 
 
\begin{lem} \label{MP_lem1} Se satisface la siguiente relacin:
\[
P_{a,h,\gamma}(i_1h_1,i_2h_2,\cdots,i_{n+1}h_{n+1})=f(a+\sum_{k=1}^{n+1}i_kh_kv_k), \text{ para  todo } (i_1,\cdots,i_{n+1})\in\mathbb{Z}^{n+1}.
\]
\end{lem}

\noindent \textbf{Demostracin. } Fijemos los valores de $k\in\{1,\cdots,n+1\}$ y $i_1$, $i_2$, $\cdots$, $i_{k-1}$, $i_{k+1}$, $\cdots$, $i_{n+1}\in\{0,1,\cdots,m\}$, y consideremos el polinomio 
$$q_k(x)=P_{a,h,\gamma}(i_1h_1,\cdots, i_{k-1}h_{k-1},x,i_{k+1}h_{k+1}, \cdots,i_{n+1}h_{n+1}) .$$ 
Es evidente que $q_k\in\Pi_m^1$ y, por tanto,
\begin{eqnarray*}
0 &=& \Delta_{h_k}^{m+1}q_k(0)=\sum_{r=0}^{m+1}\binom{m+1}{r}(-1)^{m+1-r}q_k(rh_k) \\
&=& \sum_{r=0}^{m}\binom{m+1}{r}(-1)^{m+1-r}P_{a,h,\gamma}(i_1h_1,\cdots, i_{k-1}h_{k-1},rh_k,i_{k+1}h_{k+1}, \cdots,i_{n+1}h_{n+1}) \\
&\ & \ \ \ + q_k((m+1)h_k)\\
&=& \sum_{r=0}^{m}\binom{m+1}{r}(-1)^{m+1-r}f(a+\sum_{(0\leq j\leq n+1;\ j\neq k)} i_jh_jv_j+ rh_kv_k) + q_k((m+1)h_k)\\
&=& \Delta_{h_kv_k}^{m+1}f(a+\sum_{(0\leq j\leq n+1;\ j\neq k)} i_jh_jv_j) -f(a+\sum_{(0\leq j\leq n+1;\ j\neq k)} i_jh_jv_j+ (m+1)h_kv_k)\\
&\ & \ \ \ + q_k((m+1)h_k)\\
&=&  q_k((m+1)h_k)-f(a+\sum_{(0\leq j\leq n+1;\ j\neq k)} i_jh_jv_j+ (m+1)h_kv_k).
\end{eqnarray*}
En consecuencia, tenemos que 
\begin{eqnarray*}
q_k((m+1)h_k) &=&  P_{a,h,\gamma}(i_1h_1,\cdots, i_{k-1}h_{k-1},(m+1)h_k,i_{k+1}h_{k+1}, \cdots,i_{n+1}h_{n+1}) \\
&=& f(a+\sum_{(0\leq j\leq n+1;\ j\neq k)} i_jh_jv_j+ (m+1)h_kv_k)
\end{eqnarray*}
y, repitiendo el mismo tipo de argumento (tanto hacia delante como hacia atrs, y para cada $k\in\{1,\cdots,n+1\}$), obtenemos que
\[
P_{a,h,\gamma}(i_1h_1,i_2h_2,\cdots,i_{n+1}h_{n+1})=f(a+\sum_{k=1}^{n+1}i_kh_kv_k), \text{ para  todo } (i_1,\cdots,i_{n+1})\in\mathbb{Z}^{n+1},
\] 
que es lo que queramos probar. {\hfill $\Box$}


\begin{lem}\label{MP_lem2} Se satisface la siguiente relacin:
\[
P_{a,h,\gamma}(r_1h_1,r_2h_2,\cdots,r_{n+1}h_{n+1})=f(a+\sum_{k=1}^{n+1}r_kh_kv_k), \text{ para  todo } (r_1,\cdots,r_{n+1})\in\mathbb{Q}^{n+1}.
\]
Por tanto, $\overline{G(f)}^{\mathbb{R}^{n+1}}$ contiene al conjunto $\varphi_{\gamma}(\mathbb{R}^{n+1})$, donde 
\[
\varphi_{\gamma}(t_1,\cdots,t_{n+1})=(a+\sum_{k=1}^{n+1}t_kv_k,P_{a,h,\gamma}(t_1,\cdots,t_{n+1})). 
\]
\end{lem}
\noindent \textbf{Demostracin. } Basta tener en cuenta que, si $p_1,p_2,\cdots,p_{n+1}\in \mathbb{Z}\setminus\{0\}$, y aplicamos el Lema \ref{MP_lem1}  al polinomio $P^*(t_1,\cdots,t_{n+1})$ que satisface las relaciones de interpolacin
\[
P^*(i_1h_1^*,i_2h_2^*,\cdots,i_{n+1}h_{n+1}^*)=f(a+\sum_{k=1}^{n+1}i_kh_k^*v_k), \text{ para } 0\leq i_k\leq m,\ 1\leq k\leq n+1,
\] 
donde $h_k^*=h_k/p_k$, $k=1,\cdots,n+1$, entonces 
\begin{eqnarray*}
P^*(i_1h_1,i_2h_2,\cdots,i_{n+1}h_{n+1}) &=& P^*(p_1i_1h_1^*,p_2i_2h_2^*,\cdots,p_{n+1}i_{n+1}h_{n+1}^*) \\
&=& f(a+\sum_{k=1}^{n+1}p_ki_kh_k^*v_k)\\ 
&=& f(a+\sum_{k=1}^{n+1}i_kh_kv_k)\\
&=& P_{a,h,\gamma}(i_1h_1,i_2h_2,\cdots,i_{n+1}h_{n+1}), 
\end{eqnarray*} 
para $0\leq i_k\leq m$ y $1\leq k\leq n+1$, por lo que $P^*=P_{a,h,\gamma}$, lo cual nos conduce a la identidad buscada, pues   $p_1,p_2,\cdots,p_{n+1}\in \mathbb{Z}\setminus\{0\}$ eran arbitrarios. {\hfill $\Box$}

\begin{lem}\label{MP_lem3} Si imponemos $v_k=e_k$ para $k=1,2,\cdots,n$ y $v_{n+1}=e_s$, donde $$e_i=(0,0,\cdots,1^{(\text{i-\'{e}sima posicin})},0,\cdots,0), \ \ i=1,\cdots,n,$$ 
entonces
\[
\varphi_{\gamma}(t_1,\cdots,t_{n+1})=(a+(t_1,\cdots,t_{s-1},t_s+t_{n+1},t_{s+1},\cdots,t_n),P_{a,h,\gamma}(t_1,\cdots,t_{n+1})),
\]
y
\[
P_{a,h,\gamma}(0,\cdots,0,t_s,0,\cdots,0,t_{n+1})=p_{a_s,h_s,h_{n+1}}(t_s,t_{n+1})=\sum_{k=0}^NA_k(t_s+t_{n+1})t_s^k.
\]
\end{lem}

\noindent \textbf{Demostracin. } Es trivial. Se trata simplemente de imponer las sustituciones indicadas (i.e., $v_k=e_k$ para $k=1,2,\cdots,n$ y $v_{n+1}=e_s$) y utilizar la definicin del polinomio $p_{a_s,h_s,h_{n+1}}$. {\hfill $\Box$}

El siguiente lema es un resultado conocido en Geometra Algebraica. Incluimos la demostracin con el objetivo de que nuestra prueba del Teorema \ref{MP_grafo_vv} quede establecida con nitidez en todos sus aspectos.

\begin{lem} \label{nuevoGA} Supongamos que $V\subseteq \mathbb{R}^n$ es una variedad algebraica, $V\neq \mathbb{R}^n$. Entonces $V$ no tiene puntos interiores. 
\end{lem}
\noindent \textbf{Demostracin. } Un sencillo cambio de coordenadas permite reducir la cuestin a conocer si el origen de coordenadas puede ser un punto interior de $V$. Procedemos por induccin sobre $n$. Para $n=1$ el resultado es claro (es decir, si $0\in\text{\textbf{Int}}(V)$, entonces $V=\mathbb{R}$). Supongamos que el lema es cierto para $n$ y asumamos que $\mathbf{0}=(0,\cdots,0)\in V\subset \mathbb{R}^{n+1}$. Si $\mathbf{0}\in \text{\textbf{Int}}(V)$ y $H\subseteq \mathbb{R}^{n+1}$ es un hiperplano que pasa por el origen de coordenadas, entonces   $W=V\cap H$ es una variedad algebraica que, al estar contenida en $H\simeq \mathbb{R}^n$, podemos interpretar de forma natural como una subvariedad algebraica de $\mathbb{R}^n$. Adems, el origen de coordenadas es un punto interior de $W$ (en la topologa heredada de $H$, que es la que tiene como subvariedad de $\mathbb{R}^n$). Por tanto, la hiptesis de induccin nos dice que $V\cap H=W=H$. Como esto es cierto para todos los hiperplanos $H$  que contienen a $\mathbf{0}$, y stos recubren $\mathbb{R}^{n+1}$, concluimos que $V=\mathbb{R}^{n+1}$ {\hfill $\Box$} 

\medskip

\noindent \textbf{Demostracin del Teorema \ref{MP_grafo_vv} } Utilizando la primera igualdad del Lema \ref{MP_lem3}, se sigue que
\begin{equation}\label{MP_det}
\varphi_{\gamma}'=\left [
\begin{array}{cccccccc}
1 & 0 & 0 & \cdots & 0 &\cdots & 0 & 0 \\
0 & 1 & 0 & \cdots & 0& \cdots & 0 &  0 \\
\vdots & \vdots & \vdots & \ddots & \vdots &\cdots & \vdots & \vdots \\
0 & 0 & 0 & \cdots & 1 & \cdots & 0 &  1\\
\vdots & \vdots & \vdots & \ddots & \vdots &\cdots & \vdots & \vdots \\
0 & 0 & 0 & \cdots & 0 &\cdots & 1 &  0\\
\frac{\partial P_{a,h,\gamma}}{\partial t_1} & \frac{\partial P_{a,h,\gamma}}{\partial t_2}  & \frac{\partial P_{a,h,\gamma}}{\partial t_3}  & \cdots & \frac{\partial P_{a,h,\gamma}}{\partial t_s} & \cdots & \frac{\partial P_{a,h,\gamma}}{\partial t_n}  &  \frac{\partial P_{a,h,\gamma}}{\partial t_{n+1}} 
\end{array} \right],
\end{equation}
y, por tanto, desarrollando  $\det\varphi_{\gamma}'$ por la ltima fila, obtenemos que, usando la notacin abreviada $P=P_{a,h,\gamma}$,
\begin{eqnarray*}
&\ & \xi(t_1,\cdots,t_{n+1}) := \det\varphi_{\gamma}' (t_1,\cdots,t_{n+1})\\
&\ & \ \ =  (-1)^{n+1+s}\frac{\partial P}{\partial t_{s}}(t_1,\cdots,t_{n+1})\cdot (-1)^{n-s}+\frac{\partial P}{\partial t_{n+1}}(t_1,\cdots,t_{n+1})\\
 &\ & \ \   = (\frac{\partial P}{\partial t_{n+1}} - \frac{\partial P}{\partial t_{s}})(t_1,\cdots,t_{n+1}) .
\end{eqnarray*}
Si evaluamos el polinomio $\xi$ en el punto $(0,0,\cdots,0,t_s,0,\cdots,t_{n+1})$ y usamos la segunda igualdad del Lema \ref{MP_lem3}, entonces tenemos que
\[
\det\varphi_{\gamma}'(0,0,\cdots,0,t_s,0,\cdots,t_{n+1})=- \sum_{k=1}^NkA_k(t_s+t_{n+1})t_s^{k-1}\neq 0,
\]
es un polinomio algebraico ordinario no nulo  en las variables $t_1,\cdots,t_{n+1}$. Por tanto, si tenemos en cuenta el Lema \ref{nuevoGA}, podemos afirmar que la variedad algebraica asociada,
\[
Z(\det\varphi_{\gamma}')=\{(\alpha_1,\cdots,\alpha_{n+1})\in\mathbb{R}^{n+1}: \det\varphi_{\gamma}'(\alpha_1,\cdots,\alpha_{n+1})=0\}
\]
es un subconjunto cerrado propio de $\mathbb{R}^{n+1}$ con interior vaco. Se sigue que $\Omega=\mathbb{R}^{n+1}\setminus Z(\det\varphi_{\gamma}')$ es un abierto no acotado y, gracias al Teorema de la aplicacin abierta para funciones diferenciables en dimensin finita,  
$\varphi_{\gamma}(\Omega)$ es un abierto  de $\mathbb{R}^{n+1}$ que est contenido en $\overline{G(f)}^{\mathbb{R}^{n+1}}$, que es lo que queramos probar.  La parte de no acotacin de $\varphi_{\gamma}(\Omega)$  se sigue directamente de la segunda igualdad del Lema \ref{MP_lem3}. 

{\hfill $\Box$}

\noindent \textbf{Agradecimientos. } Los autores de este trabajo agradecen al \'{a}rbitro el trabajo realizado, que ha permitido mejorar sensiblemente la redacci\'{o}n del mismo.  

 \bibliographystyle{amsplain}


\end{document}